\documentclass[12pt]{amsart}

\usepackage[utf8]{inputenc}
\usepackage{amsmath, amssymb, amsthm, mathrsfs}
\usepackage{graphicx}
\usepackage{enumitem}
\usepackage[margin=1in]{geometry}
\usepackage{natbib}
\usepackage[hidelinks]{hyperref}
\usepackage{booktabs}
\usepackage{fancyhdr}
\usepackage{setspace}
\usepackage{centernot}
\usepackage{float}

\usepackage{enumitem}

\newlist{conditions}{description}{1}
\setlist[conditions]{
  leftmargin=4.4em,
  labelwidth=3.6em,
  labelsep=0.6em,
  align=right,
  itemsep=0.35em
}
\def\s{\sigma}
\def\ve{\varepsilon}
\def\eproof{\hfill $\square$}
\def\bproof{\noindent {\em Proof. }}

\theoremstyle{definition}
\newtheorem{defi}{Definition}

\theoremstyle{plain}
\newtheorem{thm}{Theorem}
\newtheorem{lem}{Lemma}

\newtheorem{result}{Result}

\theoremstyle{remark}
\newtheorem{remark}{Remark}

\title{Robust Optimization in Stochastic Games}
\author{Abraham Neyman$^*$}
\thanks{$^*$Institute of Mathematics and the Federmann Center for the Study of Rationality, The Hebrew University of Jerusalem, Edmund Safra Campus, Jerusalem 91904, Israel. \texttt{aneyman@mail.huji.ac.il}.
\newline This research was supported in part by Israel Science Foundation grants 1596/10 and 992/23.
The author thanks Elchanan Ben-Porath, Sergiu Hart, Elon Kohlberg, Eilon Solan, Sylvain Sorin, Yair Tauman, and Bruno Ziliotto for helpful comments.}
\date{July 22, 2026}

\begin{document}

\begin{titlepage}
\maketitle

\thispagestyle{empty}

\begin{abstract}
A player's impatient valuation of an infinite stream of stage payoffs is a weighted average of those payoffs with nonincreasing weights. Impatient valuations generalize geometric discounting, and a set of them captures imprecise or heterogeneous time preferences. For a stochastic game \(\Gamma\) and a profile \(u\) of such valuations,
let \(\Gamma_u\) be the induced game. Because the valuations are rarely known exactly, we ask when approximately optimal or equilibrium behavior can be chosen to vary continuously and stably with \(u\), so that small errors in one's own or others' time preferences leave strategies and payoffs nearly unchanged.

We show that this robustness is governed by the corresponding uniform long-run payoff property: the uniform value in the adversary case, and uniform approximate equilibrium payoffs in the equilibrium case. The latter property requires that, for every tolerance, one can choose a bounded payoff vector and a single strategy profile such that, for all sufficiently long finite horizons, the profile is an approximate equilibrium and its payoffs are uniformly close to that vector. In
finite state spaces, this property is equivalent to the usual existence of a uniform equilibrium payoff.

In the two-player zero-sum case, if \(\Gamma\) has a uniform value, then the value of \(\Gamma_u\) varies uniformly continuously with \(u\), and, moreover, for every \(\varepsilon>0\), there exists a uniformly continuous mapping assigning to each \(u\) an \(\varepsilon\)-optimal strategy in \(\Gamma_u\). This holds for all
two-player zero-sum stochastic games with finite state and action sets.

In the non-zero-sum case, the existence of uniform approximate equilibrium payoffs implies that, for every \(\varepsilon>0\), the space of impatient-valuation profiles is covered by finitely many neighborhoods, on each of which a single strategy profile is a robust \(\varepsilon\)-equilibrium. 

For sunspot, correlated, and extensive-form correlated equilibria, the conclusion strengthens. For each of these classes, the existence of uniform approximate equilibrium payoffs in that class is equivalent to the existence, for every \(\varepsilon>0\), of a uniformly continuous mapping assigning to each impatient-valuation profile \(u\) a robust \(\varepsilon\)-equilibrium in that class for a neighborhood of \(u\). Consequently, every finite stochastic game admits uniformly continuous robust extensive-form correlated \(\varepsilon\)-equilibrium selections.

\end{abstract}

\end{titlepage}

\section{Introduction}
Economic models necessarily involve approximations, and time preferences are often difficult to specify precisely. In dynamic environments, small imprecisions in how an agent discounts the future---or perceives the discounting of others---can unravel classically optimal strategies. Optimizing against a single, perfectly specified model is therefore an unrealistic benchmark; a more robust approach seeks behavior that performs reliably across an entire envelope of plausible preference specifications.

This challenge amplifies when the agent is a collective entity, such as a firm, a regulatory body, or a government, whose members harbor heterogeneous time preferences. Robust optimization then requires evaluating behavior across a range of heterogeneous objectives, without assuming that the collective speaks with a single, unanimous voice.

Crucially, robust optimization also demands behavioral \emph{stability}: because shifts due to data noise, structural errors, or evolving information are unavoidable, optimizing behavior and its resulting payoffs must respond gradually---not discontinuously or erratically---to small variations in preferences.

We formalize robust optimization as behavior satisfying two core properties across both two-player zero-sum (henceforth, zero-sum) and multiplayer non-zero-sum settings:
\begin{enumerate}[label=(\roman*)]
    \item approximate optimization across a target set of preferences, accommodating local heterogeneity or modeling imprecision,
    \item gradual responsiveness (stability) of strategies and the corresponding payoffs to small shifts within that set.
\end{enumerate}

Although this criterion requires a single strategy to perform well \emph{simultaneously} across a set of preferences---a worst-case requirement over that set---it departs from the worst-case paradigm of model-misspecification analysis in two respects: it targets \emph{approximate} rather than exact optimization, and it additionally demands \emph{stability} of the induced payoffs. It likewise differs from robust implementation, which studies mechanisms achieving exact optimization robust to the information structure. Moreover, while most prior robustness studies focus on static environments, our contribution develops a robust optimization framework for dynamic, multistage interactions. We review these distinctions in detail in Section~\ref{sec:relatedliterature}.

We study robust optimization within the framework of stochastic games, a class of multiplayer, multistage interactions where actions at each stage alter both immediate payoffs and subsequent state dynamics. This intertemporal structure introduces a fundamental trade-off: agents must weigh short-term rewards against the long-term consequences of their choices. Such trade-offs govern critical economic domains, including natural resource management, dynamic bargaining, innovation investments, reputation-building, and macroeconomic policy design.

Climate policy provides a natural example. Each country chooses how much to abate over time; abatement is costly today but lowers the chance of transitions to persistently worse climate states. Because emissions are a global externality, every country's payoff---and its best response---depends on what the others do, while each country's willingness to bear present costs for future gains is governed by how it discounts the future, that is, by its impatient valuation. A country thus optimizes not against a fixed environment but against other countries whose choices are driven by their discounting. These discount rates are rarely known precisely, differ across countries, and may themselves reflect collectives---agencies, generations, coalitions---that do not speak with a single voice. Robust behavior must therefore be insensitive to imprecision both in a country's own valuation and in the profile of the others': it should remain an approximate equilibrium across the whole range of valuation profiles. This paper characterizes when such robustness is achievable.

In multistage interactions, players' preferences typically concern the entire stream of stage payoffs. We focus on preferences represented by impatient valuations, where the weight assigned to a later stage is no greater than the weight assigned to an earlier one. This class generalizes geometric discounting, requiring only that the stage weights be non-increasing.

Each impatient valuation can be interpreted as assigning a \emph{normalized present value}---a weighted average of stage payoffs with nonincreasing weights summing to one. Under this interpretation, the distance between two impatient valuations reflects the maximal discrepancy in the normalized present values assigned to any distribution of a unit total payoff spread across stages, providing an economic rationale for our metric. Familiar examples include exponential geometric discounting and time-dependent discounting where weights decline proportionally to a power $\alpha>1$ of the stage number.

Given a stochastic game $\Gamma$ and a profile $u$ of impatient valuations, $\Gamma_u$ denotes the game where each player evaluates payoffs according to their component of $u$. When $\Gamma$ is zero-sum, we write $\mathrm{val}(\Gamma_u)$ for its value and $\gamma_u(\sigma)$ for the induced payoff under a strategy $\sigma$; in general-sum settings, $\gamma_u(\sigma)$ denotes the resulting payoff profile. Unless stated otherwise, strategies are understood to be fully history-dependent behavioral strategies, which we define formally alongside memoryless variants in Section~\ref{sec:basicconcepts}.

In many economic environments, assuming that agents' impatient valuations are precisely known is unrealistic. Imprecision may arise from limited information, ambiguity about future events, or internal preference heterogeneity within a collective decision maker. To capture this uncertainty, we model preferences through a set $U$ of plausible impatient valuations or valuation profiles, requiring optimization to be robust across all elements in $U$. Such uncertainty is salient in settings like the climate policy discussed above: contested discount rates, and stakeholders who weigh short- and long-term objectives differently, translate directly into a set U of plausible valuations rather than a single one.

Two responses to such imprecision are possible. One is to aggregate the plausible valuations into a single one---for instance, by averaging geometric discount factors over a distribution of rates, which yields the declining-rate schedules advocated by Weitzman (2001) and Gollier and Weitzman (2010)\nocite{weitzman2001gamma}\nocite{gollierweitzman2010}. But aggregation need not be faithful: when \(U\) reflects a heterogeneous collective, no single time-consistent preference can aggregate its members while respecting unanimity and non-dictatorship (Jackson and Yariv, 2015)\nocite{jacksonyariv2015}.\footnote{The impossibility depends on the setting: for intergenerational social choice, Feng and Ke (2018)\nocite{fengke2018} characterize conditions under which a Paretian, non-dictatorial social planner can nonetheless discount exponentially (and hence time-consistently).} We therefore take the other route---retaining \(U\) in its entirety and seeking behavior that is robust across all of its elements.

We evaluate robust optimization in both zero-sum and non-zero-sum stochastic games. The zero-sum setting serves as an analytically tractable starting point that provides key insights extended later to general-sum games. Beyond their intrinsic interest, zero-sum dynamics are central to analyzing general multiplayer settings: determining what a player can guarantee amounts to studying the maxmin of a zero-sum game against a hypothetical adversary representing the collective of opponents. Furthermore, core challenges such as preference imprecision and the limitations of memoryless strategies transcend the zero-sum boundary, extending to Nash and correlated equilibria in general-sum games.

A strategy is \emph{\(\varepsilon\)-optimal for a set \(U\)} of impatient valuations if it is \(\varepsilon\)-optimal in \(\Gamma_u\) for every \(u \in U\). The zero-sum game \(\Gamma\) is said to have \emph{robust optimal strategies} if, for every \(\varepsilon > 0\), there is a \(\delta > 0\) such that on every set \(U\) of impatient valuations of diameter less than \(\delta\), each player~\(i\) (\(i = 1,2\)) has a single strategy \(\sigma^i\) that is \(\varepsilon\)-optimal for all \(u \in U\) and whose induced payoff \(\gamma_u(\sigma)\),
with \(\sigma = (\sigma^1,\sigma^2)\), varies by at most \(\varepsilon\) across \(U\)---that
is, \(|\gamma_u(\sigma) - \gamma_{u'}(\sigma)| < \varepsilon\) for all \(u, u' \in U\). The first requirement is approximate optimality throughout \(U\); the second is stability of the induced payoffs. Dropping it yields the weaker notion of \emph{weakly robust optimal
strategies}, which requires only that each player have an \(\varepsilon\)-optimal strategy for \(U\). The formal definitions are given in
Section~\ref{sec:robustdefs}.%

Controlling the variation of payoffs $\gamma_u(\sigma)$ within $U$ is particularly critical when $\Gamma$ arises as a subgame within a larger dynamic interaction. In these contexts, players rely on approximate backward induction, substituting the subgame with an approximate value estimating future continuation payoffs. If small changes in $u$ generate large fluctuations in payoffs, estimating future payoffs becomes impossible. Conversely, if payoffs vary only slightly across $U$, following the robust strategy $\sigma$ provides a stable proxy for continuation values.

Building on these ideas, we establish three results for finite zero-sum stochastic games (games with finitely many states and available actions):
\begin{enumerate}
    \item Any finite zero-sum stochastic game $\Gamma$ has robust optimal strategies.
    \item The value of $\Gamma_u$ is \emph{uniformly} continuous in $u$.
    \item For every $\varepsilon > 0$, there exists a uniformly continuous function mapping each impatient valuation $u$ to an $\varepsilon$-optimal strategy $\sigma_u$ in $\Gamma_u$.
\end{enumerate}

The first two results follow directly from the third, which represents a stronger form of robustness: a uniformly continuous selection mapping $u \mapsto \sigma_u$ measured under the total variation distance for mixed strategies. While standard continuity of the value mapping $u \mapsto \text{value}(\Gamma_u)$ is straightforward (see Remark~\ref{rem:L1_payoff_continuity} in Section~4), \emph{uniform} continuity requires a more delicate analysis and is indispensable for our robustness properties.

A striking feature of robust optimization in stochastic games is the limitation of Markovian strategies. For a fixed impatient valuation, finite zero-sum stochastic games admit optimal Markovian strategies. By contrast, robustness to small preference imprecision may require history-dependent behavior. The canonical example is the Big Match, where the classical uniform-optimality analysis already shows that long-run guarantees cannot, in general, be implemented by Markovian or stationary rules alone.

This phenomenon is related to the classical sunk-cost principle. That principle may be understood as requiring that, after every history, behavior solve the continuation problem from the current state, with past payoffs affecting future behavior only through their effect on the current state and primitives. Under preference imprecision, however, history can play an additional role. A strategy that is required to be approximately optimal for all valuations in a neighborhood may need to condition future play on the realized history in order to balance its performance across those valuations. When a realized path is favorable for some valuations in the neighborhood, continuation play can be adjusted so as to maintain approximate optimality for the others. Thus, robust optimization may require genuinely history-dependent strategies even in Markovian environments.

In certain structured classes of stochastic games, robust optimization can indeed be implemented via Markovian strategies. For a single decision maker (Markov decision processes), this is shown in \cite{neyman2023}; there the patient tail can be carried by a stationary strategy, so robustness lives within the memoryless class. In the game setting, by contrast, the uniform value generally requires history-dependent strategies on the tail, which is precisely what forces genuinely non-Markovian robust play: the technique of approximating by a finite horizon with a terminal continuation value is common to both, but this phenomenon is new. In finite zero-sum stochastic games with perfect information, a result from this paper combined with \cite[Theorem~1]{liggett_lippman_1969} guarantees robust optimization via Markovian strategies. Similarly, in irreducible finite zero-sum games, combining our results with \cite{HoffmanKarp66} yields the same conclusion. In the non-zero-sum setting, our framework paired with \cite{sobel1971noncooperative} establishes that robust optimization can be achieved by memoryless strategies in finite irreducible stochastic games. Outside these highly structured environments, however, genuinely history-dependent play is indispensable.

We turn to the multiplayer case. A strategy profile is an \emph{\(\varepsilon\)-equilibrium for a set \(U\)} of impatient valuation profiles if it is an \(\varepsilon\)-equilibrium of \(\Gamma_u\) for every \(u \in U\). The game \(\Gamma\) is said to have \emph{robust
equilibria} if, for every \(\varepsilon > 0\), there is a \(\delta > 0\) such that on every set \(U\) of profiles of diameter less than \(\delta\)---where profiles are compared by the maximum distance between players' component valuations---there is a single strategy profile
\(\sigma_U\) that is an \(\varepsilon\)-equilibrium of \(\Gamma_u\) for every \(u \in U\) and whose induced payoff profile varies by at most \(\varepsilon\) across \(U\). Here too, the defining requirements are approximate equilibrium throughout \(U\) together with stability of
the induced payoffs; dropping the latter yields the weaker notion of \emph{weakly robust equilibria}. The formal definitions appear in
Section~\ref{sec:robustdefs}.%

As in the zero-sum case, controlling payoff variation across $U$ enables players to utilize approximate backward induction in meta-games, treating the equilibrium payoff as a stable proxy for continuation values without introducing strategic coordination breakdowns under slight informational noise. We prove more generally that a stochastic game \(\Gamma\) has robust equilibria whenever it has uniform approximate equilibrium payoffs. This condition is weaker than the classical existence of a uniform equilibrium payoff in general, and is equivalent to it when the state space is finite. Thus the classical uniform equilibrium payoff results imply robustness in the familiar finite-state classes,
while the characterization itself is naturally formulated in terms of uniform approximate equilibrium payoffs. Several major classes of games satisfy this condition, including finite two-player stochastic games~\cite{vie00a,vie00b}, finite three-player absorbing games~\cite{solan99}, and games with perfect information~\cite{thurag97}. Whether every finite stochastic game possesses a uniform equilibrium payoff remains a foundational open problem in game theory.

Our analysis of robust optimization extends naturally to extensive-form correlated (EFC) equilibria, which generalize both correlated and sunspot equilibria by serving as external coordination mechanisms that generate incentive-compatible recommendations over time. In an EFC profile, asymmetric private signals are dispatched at each stage, allowing play to adapt dynamically to the joint history of actions and signals. An EFC equilibrium satisfies sequential incentive compatibility; an $\varepsilon$-equilibrium satisfies $\varepsilon$-incentive compatibility. Correlated equilibria restrict private signals to a single draw before play begins, while sunspot equilibria further specialize coordination to a single, publicly observed random signal revealed prior to play.

The existence of robust optimal strategies and robust equilibrium concepts (Nash, sunspot, correlated, or EFC) is characterized by the corresponding uniform long-run payoff condition. In the adversary case, robust maxmin strategies are equivalent to the existence of the uniform maxmin value. In the Nash case, robust equilibria are equivalent to the existence of uniform approximate equilibrium payoffs. In the mediated case, robust \(\mathcal M\) equilibria are equivalent to the existence of uniform approximate \(\mathcal M\) equilibrium payoffs. For finite state spaces, these approximate conditions are equivalent to the usual existence of the corresponding uniform equilibrium payoff.

While the “only if” direction follows from general principles, our main
contribution is the converse: the corresponding uniform long-run payoff
condition implies the existence of robust optimal strategies or robust
equilibria. In the adversary and mediated cases, it also yields uniformly
continuous robust selections.

These implications rely on the following topological property: the metric completion of the space of impatient valuations---the minimal extension where every Cauchy sequence converges---is compact. In the zero-sum case and for the convex equilibrium concepts (sunspot, correlated, and EFC), we exploit this compactness alongside the convexity of the approximately optimal strategy sets to construct uniformly continuous selections over the entire preference space via a partition-of-unity argument.

Broadly speaking, our results fall into two strategic categories:
\begin{itemize}
    \item \textbf{Local robustness}: If a game has uniform approximate equilibrium payoffs, we obtain local neighborhoods where an $\varepsilon$-equilibrium exists and yields stable payoffs.
    \item \textbf{Global robustness}: In zero-sum and convex non-zero-sum environments (sunspot, correlated, and EFC), we establish globally uniformly continuous selections of approximately optimal or equilibrium strategies across the entire space.
\end{itemize}

\paragraph{\textbf{Remark (On Continuity, Uniform Continuity, and Robustness).}}
The distinction between continuity and uniform continuity is central to our framework. Ordinary continuity lets the tolerance required for a given accuracy vary arbitrarily across the domain---$f(x)=\sin(1/x)$ on the open interval $(0,1)$ is continuous yet, for every $\delta>0$, has points within $\delta$ at which $|f(x)-f(x')|=2$. Simple continuity of the value or equilibrium selection therefore cannot guarantee robustness, since minor preference errors could still cause large payoff or behavioral swings. Uniform continuity rules this out, mapping uniformly small preference shifts to uniformly small outcome deviations---the natural formalization of structural stability under preference imprecision.

The paper is organized as follows. Section~\ref{sec:relatedliterature} reviews the related literature. Section~\ref{sec:basicconcepts} introduces foundational concepts. Section~\ref{sec:results} states the main results, which are highlighted in tabular form in Section~\ref{sec:summary}. Section~\ref{sec:proofofthms} details the proofs. Section~\ref{sec:othermodels} extends the analysis to other models of multistage games. Finally, Section~\ref{sec:transitionsrobustness} raises an open problem regarding robustness to imprecisions in transition probabilities.

\section{Related literature}\label{sec:relatedliterature}
\subsection{Stochastic games}
The theory of stochastic games was initiated by Shapley (1953), who introduced finite-state, finite-action two-player zero-sum stochastic games with discounted payoffs and established the existence of the value. His work inspired a broad line of research extending to stochastic games with long-run average payoffs.

Two main approaches have guided the study of long-run average payoffs in stochastic games: the \emph{undiscounted approach} and the \emph{uniform approach}. Under the undiscounted approach, the objective is to optimize a specific long-run average payoff criterion, such as the limit inferior or limit superior of the average payoff over the first $n$ stages. By contrast, the uniform approach considers performance across all sufficiently long finite-horizon games or sufficiently patient discounted games, aiming for strategies that guarantee nearly optimal outcomes in every such environment. Historically, much of the early work---including that of Hoffman and Karp (1966), Liggett and Lippman (1969), and Blackwell and Ferguson (1968)---focused on the existence of the value under the undiscounted approach.
Although these results were formulated for lim-average payoffs, their proofs construct strategies with uniform $\varepsilon$-optimality properties: for every $\varepsilon>0$, the strategies can be chosen to be $\varepsilon$-optimal uniformly over all sufficiently long horizons.

Early results within the undiscounted approach established the existence of the value for several special classes of stochastic games. For example, Hoffman and Karp (1966)\nocite{HoffmanKarp66} proved the existence of the value and optimal stationary strategies in finite irreducible stochastic games, and Liggett and Lippman (1969)\nocite{liggett_lippman_1969} established the same result for finite stochastic games with perfect information.

Gillette (1957)\nocite{gil57} introduced a simple stochastic game—later termed the Big Match—that revealed a striking feature of stochastic games: in some cases, players must use complex strategies that condition their behavior on the full history of play in order to perform near-optimally. The Big Match became a classic example illustrating the limitations of memoryless or stationary strategies. It was subsequently solved by Blackwell and Ferguson (1968),\nocite{blafer68} who established the existence of the value under the limsup-average payoff.

Subsequently, Kohlberg (1974)\nocite{koh74} proved the existence of the value for the broader class of absorbing games, which generalize games like the Big Match.

The above-mentioned results on the existence of the undiscounted value highlighted the structural complexity of stochastic games and the challenges of extending existence results beyond highly specialized settings.

Bewley and Kohlberg (1976)\nocite{bewkoh76} introduced an algebraic approach to the analysis of finite stochastic games and used it to establish the existence and equality of two limits: the limit of the $n$-stage values as $n \to \infty$, and the limit of the $\lambda$-discounted values as $\lambda \to 0$. Their framework provided useful tools for subsequent developments.

\cite{merney81,MertensNeyman1982}
established the existence of the uniform value in finite two-player
zero-sum stochastic games. Moreover, for every $\varepsilon>0$, the same
strategies can be chosen to be $\varepsilon$-optimal uniformly for all
sufficiently long finite-horizon games, for all sufficiently patient
discounted games, and for the liminf-average payoff.

In the non-zero-sum setting, Fink (1964)\nocite{Fink1964} and Takahashi (1964)\nocite{Takahashi1964} independently extended Shapley’s results by proving the existence of Nash equilibria in discounted stochastic games under general compactness and continuity assumptions.

The existence of a uniform equilibrium payoff in finite non-zero-sum stochastic games is known only in special cases; in particular, whether it always exists in finite non-zero-sum stochastic games remains unresolved.

Significant progress on this question has been made in recent decades. In a major breakthrough in the field, Vieille (2000a, 2000b) proved the existence of uniform equilibrium payoffs in all finite two-player non-zero-sum stochastic games. Solan (1999)\nocite{solan99} established the existence of uniform equilibrium payoffs in three-player absorbing games, and Solan and Vieille (2002)\nocite{solvie02} established the existence of uniform extensive-form correlated equilibrium payoffs in finite non-zero-sum stochastic games. Nevertheless, the existence of uniform equilibrium payoffs in finite
stochastic games with three or more players remains an open and
challenging problem.

Surveys and broader expositions of stochastic games---including discussions of uniform and correlated equilibria---can be found in Solan’s \emph{A Course in Stochastic Game Theory} (2022)\nocite{solan2022course}, the edited volume \emph{Stochastic Games and Applications} by Neyman and Sorin (2003)\nocite{neyman_sorin_2003}, and the overview by Solan and Vieille (2015)\nocite{solan2015stochastic}.

The notions of uniform value and uniform equilibrium payoff already capture a form of robustness to players' time preferences. In the zero-sum case, a uniform value means that, for every \(\varepsilon>0\), a strategy can be chosen that is \(\varepsilon\)-optimal for all sufficiently patient valuations, including both long finite-horizon averages and sufficiently patient discounted payoffs. In the non-zero-sum case, a uniform equilibrium payoff similarly requires \(\varepsilon\)-equilibrium strategy profiles whose induced payoffs vary by at most \(\varepsilon\) across sufficiently patient valuations.
The set of impatient valuations generated by geometric discounting with discount rate at most \(\lambda\) and by \(n\)-stage averaging with \(n \geq 1/\lambda\) has diameter at most \(\lambda\) in the supremum \((\ell_\infty)\) norm. Hence, the existence of a uniform value or uniform equilibrium payoff guarantees robustness over these standard
sets of sufficiently patient valuations, while our results extend this guarantee to arbitrary sets of impatient valuations with sufficiently small diameter.

This paper formalizes these ideas by introducing precise definitions of robust optimization and robust equilibrium under uncertainty about time preferences. It then shows how the corresponding uniform long-run payoff condition gives rise to robustness and stability of strategic behavior across nearby impatient valuations.

\subsection{Robust Implementation and Robust Optimization}
Robustness considerations in economic theory have been extensively studied, particularly in the contexts of mechanism design and implementation theory. Early work, such as \cite{BergemannMorris2005,BergemannMorris2011}, investigates the robust implementation of mechanisms under uncertainty about agents’ private information, types, or beliefs. These mechanisms aim to ensure that a target outcome, or a desired set of outcomes, is implemented across a wide range of informational environments, provided that agents behave optimally given their private information.

Another strand of research, exemplified by Wilson’s (2006)\nocite{Wilson2006} doctrine, advocates progressively weakening the informational assumptions underlying economic models. Wilson emphasizes the need for mechanisms and strategies that remain effective under broader and less structured forms of uncertainty.

The present work is related in spirit to the robust mechanism design frameworks discussed above, but departs from them in key respects. First, we consider environments in which \emph{approximate optimization} is necessary: when there is uncertainty about agents' time preferences, exact optimization across all possible preference specifications is typically unattainable. The objective thus shifts to achieving near-optimal outcomes—uniformly across a class of plausible preference specifications—rather than exact optimality. Second, we impose a \emph{continuity requirement}: small perturbations in preferences must lead to only small changes in agents’ approximate optimizing behavior. This stability condition ensures that strategic behavior is not overly sensitive to small modeling errors.

Moreover, while most prior robustness studies have focused on static models, the present work develops a robustness framework for dynamic, multistage interactions. To operationalize this approach, we study robust optimization and equilibrium behavior within the framework of stochastic games, a canonical class of dynamic models that captures both intertemporal uncertainty and strategic interaction.

The term robust optimization is also used in operations research for a distinct, worst-case paradigm in which one optimizes against the least favorable realization of the data within an uncertainty set (Ben-Tal, El Ghaoui, and Nemirovski, 2009\nocite{bental2009robust}; Bertsimas, Brown, and Caramanis, 2011\nocite{bertsimas2011theory}). Our notion differs: rather than worst-case protection, we require approximate optimality simultaneously across a small set of impatient valuations together with stability of the induced payoffs.

\subsection{Misspecification and Robust Optimization}
Misspecification has long been a central concern in economic theory, particularly in models where agents may not fully understand their environment. In robust control theory (e.g., Hansen and Sargent, 2001; 2007),\nocite{HansenSargent2001,hansen2007robustness} agents guard against model misspecification by selecting strategies that perform well across a set of possible environments, typically using a worst-case or minimax approach. Gilboa and Schmeidler’s (1989)\nocite{gilboa1989maxmin} maxmin expected utility and Epstein and Schneider’s (2003)\nocite{EpsteinSchneider2003} dynamic extension reflect similar concerns about imprecise beliefs.

In strategic settings, related ideas appear in equilibrium concepts such as self-confirming equilibrium (Fudenberg and Levine, 2006)\nocite{FudenbergLevine2006} and learning with misspecified models (Esponda and Pouzo, 2016; Esponda, 2016).\nocite{EspondaPouzo2016,Esponda2016}

Our approach addresses a distinct form of uncertainty: uncertainty about agents' time preferences. Rather than optimizing relative to a fixed impatient valuation, players seek strategies that are approximately optimal across a small neighborhood of plausible valuations. The setting can thus be read as a structured form of misspecification---of the discount weights rather than of the environment or beliefs---but with a different objective: approximate optimality and stability across the neighborhood rather than worst-case protection.

\section{Basic concepts}\label{sec:basicconcepts}

\subsection{Valuations}
\nocite{chakraborty2021present}
In dynamic environments, agents often evaluate infinite sequences of payoffs using a single summary statistic. Such a statistic is called a \emph{valuation} if it is a real-valued linear functional \(u\), defined on the space of bounded payoff streams, that satisfies the following two conditions:
\begin{itemize}
    \item \textbf{Normalization:} \( u(\mathbf{1}) = 1 \), where \( \mathbf{1} \) is the constant stream with value 1 in every period.
    \item \textbf{Time Value of Money:} If \( g \) weakly dominates \( f \) in cumulative payoffs at every stage---that is, \( \sum_{t=1}^s g_t \geq \sum_{t=1}^s f_t \) for all \( s \geq 1 \)---then \( u(g) \geq u(f) \).
\end{itemize}

The Normalization condition serves to simplify statements of results by anchoring all constant payoff streams at a common reference level. The Time Value of Money principle, on the other hand, captures a monotonicity property: payoff streams that yield greater cumulative rewards at every point in time are assigned at least as much value. Both conditions are satisfied by common evaluation criteria such as discounting and finite averaging.

We are particularly interested in \emph{impatient valuations}, which place diminishing weight on future payoffs. A valuation \( u \) is said to be \emph{impatient} if \( u(\mathbf{1}_{>t}) \to 0 \) as \( t \to \infty \), where \( \mathbf{1}_{>t} \) denotes the stream with payoffs equal to 1 after period \(t\) and 0 in all other periods.

Every impatient valuation corresponds to a nonincreasing sequence of nonnegative weights \( w = (w_1, w_2, \ldots) \) with total mass \( \sum_{t=1}^\infty w_t = 1 \), such that
\[
u(g) = \sum_{t=1}^\infty w_t g_t.
\]
We identify the set \( W \) of all impatient valuations with the set of such normalized nonincreasing weight sequences, and denote by \( w \) both the sequence and the associated linear functional. We endow the set \( W \) with the sup-norm metric,
\[
d(w, w') = \sup_{t \geq 1} |w_t - w'_t| = \max_{t \geq 1} |w_t - w'_t|,
\]
which measures the maximum pointwise discrepancy between two weight sequences.

Two classic examples of impatient valuations are the \emph{finite-horizon average} and the \emph{discounted} valuation:
\begin{itemize}
    \item The \( n \)-stage average (Cesàro average), denoted by \( u_n \), assigns equal weight \( 1/n \) to the first \( n \) periods:
    \[
    u_n(g) = \frac{1}{n} \sum_{t=1}^n g_t.
    \]
    \item The \( r \)-discounted valuation (for discount rate \( 0 < r <1 \) and discount factor \( 0 < 1-r <1 \)), denoted by \( u_r \), assigns exponentially decreasing weights over time:
    \[
    u_r(g) = \sum_{t=1}^\infty r(1-r)^{t-1} g_t.
    \]
\end{itemize}

More generally, any impatient valuation \( w \) can be interpreted as a form of time-dependent discounting, where the local discount factor at time \( t \) is given by \( w_{t+1}/w_t \).

An average of \(r\)-discounted valuations, \(\int u_r\,d\mu(r)\) where \(\mu\) is a probability measure on \((0,1)\), is an impatient valuation with local discount factors that are nondecreasing in \(t\). Since an impatient valuation need not have non-decreasing local discount factors, not every impatient valuation is such an average. Every impatient valuation is, however, an average of the \(n\)-stage average valuations \(u_n\): with \(\lambda_n=n(w_n-w_{n+1})\geq 0\) one has \(\sum_n\lambda_n=1\) and \(w=\sum_n\lambda_n\,u_n\). Equivalently, the \(u_n\) are the extreme points of the set
of impatient valuations (see, e.g., \cite{sorin2002firstcourse}).

Other examples include valuation based on power-law decay, such as \( w_t = c_\alpha / t^\alpha \) for some \( \alpha > 1 \) (with \( c_\alpha \) chosen to normalize the sum to 1), and mixtures of geometric discounting, which arise in models with heterogeneous discounting or uncertainty over agents’ time preferences. These examples include a few commonly used impatient valuations across both classical and behavioral economic models.\footnote{Some related works refer to impatient valuations in specific forms or passing discussions; see, for example, \cite{chakraborty2021present,kimli2022folk,renaultvenel2016longterm,venelziliotto2021history}.}

\paragraph{\textbf{Notation.}} We use \(u\) and \(w\) to denote impatient valuations; each may stand, according to context, for a single valuation (an element of \(W\)) or for a profile (an element of \(W^I\)), the meaning being clear from the setting. When two valuations appear together they are distinguished as \(u\) and \(w\); in the zero-sum case the players share a common valuation, denoted \(w\). We also write \(w\) for points of the completions \(\overline{W}\) and
\(\overline{W}^I\).

\subsection{Finite stochastic games}
A finite stochastic game proceeds in stages. At each stage, the game is in one of finitely many states. Players simultaneously choose actions from finite sets, with the resulting action profile and current state jointly determining the stage payoffs and the distribution over the next state.

A finite stochastic game \( \Gamma \) is defined by the 5-tuple \( \langle I, S, A, r, q \rangle \), where:
\begin{itemize}
    \item \( I \) is the finite set of players.
    \item \( S \) is the finite set of states.
    \item \( A(z) = \times_{i \in I} A^i(z) \) is the set of action profiles available at state \( z \in S \), where \( A^i(z) \) is the set of actions available to player~\( i \) at \( z \). The set of feasible state-action pairs is denoted by \( SA := \{(z, a): z \in S,\, a \in A(z)\} \).
    \item \( r: SA \to \mathbb{R}^I \) is the stage payoff function.
    \item \( q: SA \to \Delta(S) \) is the transition function, where \( q(z' \mid z, a) \) denotes the probability that the next state is \( z' \) given current state \( z \) and action profile \( a \).
\end{itemize}

A state \( z \in S \) is \emph{absorbing} if \( q(z \mid z, a) = 1 \) for all \( a \in A(z) \); that is, the game never leaves the state \( z \), regardless of the players' actions. An \emph{absorbing game} is a stochastic game with at most one nonabsorbing state.

A one-player stochastic game is known as a \emph{Markov decision process} (MDP). A two-player zero-sum stochastic game is one in which \( I = \{1, 2\} \) and \( r^2 = -r^1 \).

A finite stochastic game \(\Gamma\) has \emph{perfect information} if, at
each state \(z\), at most one player has more than one available action;
that is,
\[
\left|\{\, i : |A^i(z)|>1 \,\}\right| \leq 1
\quad \text{for every } z\in S .
\]
Equivalently, at each state, all but possibly one player have singleton
action sets.

An \emph{irreducible stochastic game} is one for which there exists a
nonempty subset \(C\subseteq S\) such that, under every profile of pure
stationary strategies, \(C\) is the unique closed communicating class of
the induced Markov chain on \(S\).

A \emph{play} of \( \Gamma \) is a sequence \( (z_1, a_1; z_2, a_2; \ldots) \in (SA)^\mathbb{N} \), where \( z_t \) is the state and \( a_t \) the action profile at stage \( t \). The initial segment \( (z_1, a_1; \ldots; z_t) \) is called the \emph{history} up to stage \( t \).

A \emph{pure strategy} for player \(i\) selects an action \( a^i_t \in A^i(z_t) \) at each stage \( t \), conditional on the history \( (z_1, a_1; \ldots; z_t) \) up to stage \(t\). A \emph{mixed strategy} for player \(i\) is a probability distribution over his pure strategies. The space of pure strategies can be identified with a countable product of finite sets, and is naturally endowed with the product topology and its associated Borel \(\sigma\)-algebra.

A \emph{behavioral strategy} for player~\(i\) specifies, at each stage \(t\), a probability distribution over \( A^i(z_t) \), conditional on the history up to stage \(t\).
Note that a pure strategy is a special case of a behavioral strategy in which the selected distributions are Dirac measures (i.e., concentrated on a single action) at every stage.
A \emph{Markov} strategy is a behavioral strategy in which the choice of mixed action at stage~\(t\) depends only on the pair~\((t,z_t)\), and a  \emph{stationary} strategy is one in which the choice of mixed action at stage~\(t\) depends only on the current state~\(z_t\).

A \emph{strategy profile} is a tuple \( \sigma = (\sigma^i)_{i \in I} \), where each \( \sigma^i \) is a behavioral strategy for player \( i \). Let \( \Sigma^i \) denote the set of behavioral strategies for player \( i \), and define \( \Sigma := \times_{i \in I} \Sigma^i \), and \( \Sigma^{-i} := \times_{j \ne i} \Sigma^j \). Given \( \tau^i \in \Sigma^i \) and \( \sigma^{-i} \in \Sigma^{-i} \), the strategy profile in which player \( i \) uses \( \tau^i \) and each other player uses \( \sigma^j \) is denoted by \( (\tau^i, \sigma^{-i}) \).

A strategy profile \( \sigma \) and initial state \( z_1 \in S \) induce a probability distribution \( P^{z_1}_\sigma \) over plays \( (z_t, a_t)_{t \geq 1} \), defined recursively as follows:
\begin{itemize}
    \item At each stage \( t \), every player \( i \) selects a mixed action \( x^i_t \in \Delta(A^i(z_t)) \) according to \( \sigma^i \), based on the history up to stage \(t\).
    \item Then, each action \( a^i_t \) is drawn independently according to the distribution \( x^i_t \).
    \item Then, the next state \( z_{t+1} \) is drawn from the distribution \( q(\cdot \mid z_t, a_t) \).
\end{itemize}

Expectations with respect to \( P^{z_1}_\sigma \) are denoted by \( E^{z_1}_\sigma \).

Because stochastic games exhibit perfect recall, for every mixed strategy there exists a%
behavioral strategy that induces the same distribution over plays, and for every behavioral strategy there exists a mixed strategy that induces the same distribution over plays.

\smallskip
\paragraph{\textbf{Metric convention for strategy spaces.}}
When discussing continuity of strategy selections, we view a mixed strategy as a probability measure over the measurable space of pure strategies, which in the finite case is the compact product space described above, and endow the space of mixed strategies with the total variation distance. A convex combination of mixed strategies is interpreted as an initial randomization over complete contingent plans. For strategy profiles, the distance is the maximum over the players' componentwise distances.

This convention is useful in the sequel because the sets of approximately maxmin mixed strategies are convex, and because total variation is a strong behavioral metric: if two mixed strategies of a player are close in total variation, then, against any fixed strategies of the other players and from any initial state, the induced distributions over plays are close in total variation. Since the games have perfect recall, mixed and behavioral strategies induce the same distributions over plays; we shall use behavioral-strategy notation when convenient, while the continuity statements are understood with respect to the total variation metric on mixed strategies.

For correlated, sunspot, and extensive-form correlated strategies, we use the analogous convention. A mediated strategy is viewed as a probability distribution over the relevant complete mediation device, including the realization of the public or private signals and the induced recommendations. A convex combination \(\sum_j \alpha_j \sigma_j\) is interpreted as follows: before play begins, an auxiliary label \(j\) is drawn with probabilities \((\alpha_j)_j\), and conditional on \(j\), the mediated strategy \(\sigma_j\) is used. Thus the label \(j\) is part of the extended randomization device. Under this convention,
\[
d\left(\sum_j \alpha_j \sigma_j,\sum_j \beta_j \sigma_j\right)
\leq
\|(\alpha_j)_j-(\beta_j)_j\|_{\mathrm{TV}},
\]
and equality holds when the common components \(\sigma_j\) are represented on disjoint labelled copies. This is the only continuity property of the metric on mediated strategies used below.

Let \( g_t = r(z_t, a_t) \) be the stage-\( t \) payoff vector, and let \( g^i_t \) denote the payoff to player \( i \) at stage \( t \).

Fix a profile of impatient valuations \( w = (w^i)_{i \in I} \in W^I \). Then for every strategy profile \( \sigma \) and initial state \( z \in S\), the \( w^i \)-evaluated payoff to player \( i \) is given by
\[
\gamma^i_w(\sigma)[z] := E^z_\sigma \left[ \sum_{t=1}^\infty w^i_t g^i_t \right].
\]
This defines the stochastic game \( \Gamma^z_w \), induced by the valuation profile \( w \) and initial state \( z \).

For the \(n\)-stage average valuation \(u_n\), we use the shorthand
\(\gamma^i_n:= \gamma^i_{u_n}\).

\subsection{Stochastic games with infinitely many states and actions}
Our framework and results extend beyond finite stochastic games, but also to more general models in which the sets of states and actions may be infinite. Specifically, we consider stochastic games with a finite set of players \( I \), a measurable state space \( S \), and a measurable action space \( A^i \) for each player \( i \in I \).

We assume the stage payoff function \( g: S \times A \to \mathbb{R}^I \) is measurable and bounded: \( \sup_{(z, a)} \|g(z, a)\| < \infty \). In addition, we require that for every bounded measurable function \( \alpha: S \to \mathbb{R}^I \), every state \( z \in S \), and every \( \varepsilon > 0 \), the one-shot game with payoff function
\[
a \mapsto g(z, a) + \int_S \alpha(z')\, q(dz' \mid z, a)
\]
admits an \( \varepsilon \)-equilibrium.

Finally, we note that even when the (finite) set of feasible states and actions vary by stage---for instance, in models where time is explicitly embedded in the  state space---the resulting game can be reformulated as a stochastic game with (countable) infinite state and action spaces. Such formulations fall within the general model described above.
Accordingly, although some robustness results rely on finiteness assumptions, the broader implication that the corresponding uniform long-run payoff condition ensures robust optimization extends to this more general class of stochastic games.

\noindent \textbf{Measurability convention.} In finite state spaces the following qualification is automatic. Whenever \(S\) is infinite, all payoff vectors \(v=(v^i[z])_{i\in I,z\in S}\) appearing below are understood to be bounded measurable functions from \(S\) to \(\mathbb R^I\), and all strategies, mediated strategy profiles, and unilateral deviations are understood to be measurable. We shall not repeat this qualification in each definition.

\subsection{The maxmin value}
We next define the maxmin, the value that a player can guarantee in a stochastic game regardless of the behavior of the other players.

\begin{defi}
The \emph{maxmin value of player \( i \)} in a stochastic game \( \Gamma_u \), where \( u \) is a profile of impatient valuations, is the vector
\[
v^i(u) = \big( v^i(u)[z] \big)_{z \in S}
\]
defined for each \( z \in S \) by
\[
v^i(u)[z] := \sup_{\sigma^i} \inf_{\sigma^{-i}} \gamma^i_u(\sigma^i, \sigma^{-i})[z],
\]
where the supremum is over all strategies of player \( i \), and the infimum is over all strategy profiles of the other players.
\end{defi}

A strategy \( \sigma^i \) is an \emph{\( \varepsilon \)-maxmin strategy} for player \( i \) in \( \Gamma_u \) if
\[
\inf_{\sigma^{-i}} \gamma^i_u(\sigma^i, \sigma^{-i})[z] \geq v^i(u)[z] - \varepsilon
\]
for every initial state \( z \in S \).

\smallskip
Note that both \( \gamma^i_u \) and \( v^i(u) \) depend only on the valuation \( u^i \) of player~\( i \), and are independent of the valuations of the other players. Accordingly, for any impatient valuation \( u^i \) of player \( i \), we may write \( \gamma^i_{u^i}\)              and \( v^i(u^i)\)  to denote the payoff function and maxmin value, respectively, corresponding to any profile \( w \) with \( w^i = u^i \).

The maxmin value of player \( i \) in the \( n \)-stage game is \( v^i(w) \) when \( w^i = u_n \), and is denoted by \( v^i_n \). In the two-player zero-sum case, we assume players share a common valuation \( w^1 = w^2=w \) and denote the maxmin value of Player~1 by \( v(w) \).

The \emph{minmax value} of player \( i \) in \( \Gamma_u \) is defined analogously:
\[
\inf_{\sigma^{-i}} \sup_{\sigma^i} \gamma^i_u(\sigma^i, \sigma^{-i})[z].
\]

In a finite stochastic game, the space of strategy profiles---viewed as a countable product of finite-dimensional simplices---is compact under the product topology, and the function \( \sigma \mapsto \gamma^i_u(\sigma) \) is continuous with respect to this topology. Hence, the above suprema and infima can be replaced by maxima and minima.

In a two-player zero-sum finite stochastic game, the maxmin and minmax values coincide for each initial state and are referred to simply as the \emph{value} of the game.

\subsection{Robust optimal strategies and equilibria of a stochastic game}\label{sec:robustdefs}
Recall that \( W \) denotes the metric space of all nonincreasing sequences \( w = (w_1, w_2, \ldots) \) with \( \sum_{t=1}^\infty w_t = 1 \), which we identify with the set of all impatient valuations.
Fix a stochastic game \( \Gamma \). Recall that the set of strategies of player \( i \) is denoted by \( \Sigma^i \), and the set of strategy profiles of the other players by \( \Sigma^{-i} \).

To model imprecision in the specification of impatient valuations we consider a set \( U \subset W \) of impatient valuations.

\begin{defi}\label{def:maxmin_strategy_for_U}
Let \( U \subset W \). A strategy \( \sigma^i \) of player~\( i \) is an \emph{\( \varepsilon \)-maxmin strategy for \( U \)} if it is an \( \varepsilon \)-maxmin strategy for every valuation \( u \in U \) of player \( i \); that is,
\begin{equation}\label{eqn:epsilonoptimalU1}
\gamma^i_u(\sigma^i, \sigma^{-i})[z] \geq v^i(u)[z] - \varepsilon
\quad \text{for all } u \in U,\, \sigma^{-i} \in \Sigma^{-i},\, z \in S.
\end{equation}
\end{defi}

Our central notion of robust optimization under imprecision in time preference specification rests on two requirements. The first (preliminary) condition  is that the strategy be an $\varepsilon$-maxmin strategy for the set \(U\); that is, the strategy must be nearly optimal for each \(u\in U\). The second, which captures the core robustness criterion and implies the first, requires that the payoff the strategy guarantees under any \(u\in U\) remains close to the maxmin payoff it guarantees under any other \(w\in U\). This motivates the following definition.

\begin{defi}\label{def:robust_maxmin_strategy_for_U}
Let \( U \subset W \). A strategy \( \sigma^i \) is a \emph{robust \( \varepsilon \)-maxmin strategy} for \( U \) if
\begin{equation}\label{eqn:epsilonoptimalforU}
\gamma^i_u(\sigma^i, \sigma^{-i})[z] \geq v^i(w)[z] - \varepsilon
\quad \text{for all } u, w \in U,\, \sigma^{-i} \in \Sigma^{-i},\, z \in S.
\end{equation}
\end{defi}

The definition yields an immediate implication regarding the variation of the maxmin value across \(U\).

\smallskip \noindent\textbf{Value variation under robust maxmin strategies.}
If \(\sigma^i\) is a robust \(\varepsilon\)-maxmin strategy for \(U\subseteq W\), then, for every \(u,w\in U\) and every \(z\in S\),
\[
v^i(u)[z]\ge v^i(w)[z]-\varepsilon .
\]
Indeed,
\[
v^i(u)[z]
\ge
\inf_{\tau\in\Sigma^{-i}}\gamma^i_u(\sigma^i,\tau)[z]
\ge
v^i(w)[z]-\varepsilon .
\]
Interchanging \(u\) and \(w\), it follows that
\[
\|v^i(u)-v^i(w)\|\le \varepsilon
\quad\text{for all }u,w\in U.
\]

The robustness results below will be derived from the following long-run value notion.

\begin{defi}\label{def:uniform_maxmin} A \emph{uniform maxmin value} of player~\(i\) is a vector \(v^i=(v^i[z])_{z\in S}\) such that for every \(\varepsilon>0\), there exist an integer \(n_\varepsilon\) and a strategy \(\sigma^i_\varepsilon\in\Sigma^i\) such that the following two conditions hold.

First, for every integer \(n\ge n_\varepsilon\), \[ \gamma_n^i(\sigma^i_\varepsilon,\tau)[z]\ge v^i[z]-\varepsilon \quad \text{for all }\tau\in\Sigma^{-i},\ z\in S. \]
Second, for every strategy \(\rho^i\in\Sigma^i\) and every impatient valuation \(w\in W\) with \(w_1<1/n_\varepsilon\), there exists \(\tau=\tau(\rho^i,w)\in\Sigma^{-i}\) such that \[ \gamma_w^i(\rho^i,\tau(\rho^i,w))[z]\le v^i[z]+\varepsilon \quad \text{for all }z\in S. \] \end{defi}

Intuitively, this means that, for any \(\varepsilon>0\), player~\(i\) has a strategy that guarantees at least \(v^i[z]-\varepsilon\), uniformly over all sufficiently long Cesaro horizons and against all opponents' strategies. Conversely, for any sufficiently patient valuation \(w\) and any strategy of player~\(i\), the other players can choose a counter-strategy, possibly depending on both \(w\) and that strategy, which keeps player~\(i\)'s payoff below \(v^i[z]+\varepsilon\).

\begin{remark}
In stochastic games with infinitely many states or actions, a uniform maxmin value may fail to exist. The definition implies that whenever a uniform maxmin value exists, it is unique.
\end{remark}

\begin{remark} By Lemma~\ref{lem:cesaro_patient}, the first condition in the definition is equivalent to the following patient-valuation formulation: for every \(\varepsilon>0\), there exist \(\delta>0\) and a strategy \(\sigma^i_\varepsilon\in\Sigma^i\) such that
\[ \gamma_w^i(\sigma^i_\varepsilon,\tau)[z]\ge v^i[z]-\varepsilon \quad \text{for all impatient valuations }w\in W\text{ with }w_1<\delta, \ \tau\in\Sigma^{-i},\ z\in S. \]
\end{remark}

In every finite stochastic game, the uniform maxmin value of each player exists in the sense of Definition~\ref{def:uniform_maxmin}; see Neyman~[2003]. The proof builds on the uniform-value theorem of Mertens and Neyman~[1981].\footnote{In fact, in finite stochastic games the proof in Neyman~[2003] yields the stronger property that, for each strategy \(\rho^i\) of player~\(i\), the opponents' counter-strategy in the second condition can be chosen independently of the sufficiently patient valuation \(w\).}

In a two-player zero-sum stochastic game, the \emph{uniform value}, when it exists, coincides with the uniform maxmin value of Player~1.

We now extend this robustness framework from individual optimization to strategic interactions, focusing on equilibrium behavior. As with maxmin strategies, our robustness criterion for equilibrium behavior rests on two requirements.

The first requirement is that the strategy profile constitutes an \(\varepsilon\)-equilibrium in every game \(\Gamma_u\) with \(u\in U\); that is, each player's strategy is \(\varepsilon\)-optimal given the others’ strategies, for every  profile of impatient valuations \(U\).

\begin{defi}
Let \( U \subset W^I \). A strategy profile \( \sigma \) is an \emph{\( \varepsilon \)-equilibrium for \( U \)} if, for every profile \( u \in U \), the strategy profile \( \sigma \) is an \( \varepsilon \)-equilibrium in the game \( \Gamma_u \); that is,
\begin{equation}\label{SGepsilonequilibriumatw}
\gamma^i_u(\sigma^{-i},\tau^i)[z] \leq \gamma^i_u(\sigma)[z] + \varepsilon
\quad \text{for all } i \in I,\, u \in U,\, \tau^i \in \Sigma^i,\, z \in S.
\end{equation}
\end{defi}

The second requirement motivates the definition of a robust \(\varepsilon\)-equilibrium for a set \(U\) of impatient valuation profiles. It requires that no unilateral deviation by a player yields a payoff---under any profile \(u\in U\)---that exceeds the payoff yielded by the strategy profile under any other profile \(w\in U\) by more than \(\varepsilon\).

\begin{defi}\label{def:robust_equilibrium_for_U}
Let \( U \subset W^I \). A strategy profile \( \sigma \) is a \emph{robust \( \varepsilon \)-equilibrium for \( U \)} if
\begin{equation}\label{SGrobustepsilonequilibriumatw}
\gamma^i_u(\sigma^{-i},\tau^i)[z] \leq \gamma^i_w(\sigma)[z] + \varepsilon
\quad \text{for all } i \in I,\, u, w \in U,\, \tau^i \in \Sigma^i,\, z \in S.
\end{equation}
\end{defi}

The definition implies that, for each player~\(i\) and initial state~\(z\), the payoff \( \gamma^i_u(\sigma)[z] \) varies by at most \( \varepsilon \) across all \( u \in U \); that is, the oscillation of \( \gamma^i_u(\sigma)[z] \) over \( u \in U \) is bounded by \( \varepsilon \).

In contrast to the maxmin value, which depends only on player \(i\)'s own valuation, an \(\varepsilon\)-equilibrium and a robust \( \varepsilon \)-equilibrium for \(U\) depend on the impatient valuations of all players.

\medskip
We now define properties of a stochastic game that captures the existence of \(\varepsilon\)-maxmin and robust \(\varepsilon\)-maxmin strategies, uniformly over all sufficiently small sets of impatient valuations.

A small imprecision in the specification of impatient valuations is modeled as a subset \( U \subset W \) with small diameter.

\begin{defi}
A stochastic game \emph{has weakly robust maxmin strategies} for player~\( i \) if, for every \( \varepsilon > 0 \), there exists \( \delta > 0 \) such that for every set \( U \subset W \) with \( \operatorname{diam}(U) < \delta \), player~\( i \) has an \( \varepsilon \)-maxmin strategy for \( U \).
\end{defi}

\begin{defi}
A stochastic game \emph{has robust maxmin strategies} for player \( i \) if, for every \( \varepsilon > 0 \), there exists \( \delta > 0 \) such that for every set \( U \subset W \) with \( \operatorname{diam}(U) < \delta \), player \( i \) has a robust \( \varepsilon \)-maxmin strategy for \( U \).
\end{defi}

\medskip
Similarly, we define properties of a stochastic game that correspond to the existence of \(\varepsilon\)-equilibria and robust \(\varepsilon\)-equilibria, uniformly over all sufficiently small sets of impatient valuation profiles.

\begin{defi}
A stochastic game \emph{has weakly robust equilibria} if, for every \( \varepsilon > 0 \), there exists \( \delta > 0 \) such that for every set \( U \subset W^I \) with \( \operatorname{diam}(U) < \delta \), there exists an \( \varepsilon \)-equilibrium for \( U \).
\end{defi}

\begin{defi}
A stochastic game \emph{has robust equilibria} if, for every \( \varepsilon > 0 \), there exists \( \delta > 0 \) such that for every set \( U \subset W^I \) with \( \operatorname{diam}(U) < \delta \), there exists a robust \( \varepsilon \)-equilibrium for \( U \).
\end{defi}

\medskip
A \emph{dummy player} is one whose set of strategies is a singleton, that is, the player has no decision to make. Similarly, in an MDP, a \emph{dummy decision maker} faces no choice. In such cases, weakly robust maxmin strategies and weakly robust equilibria trivially exist.

However, robust solutions need not exist. For example, consider an MDP with countably many states and a dummy decision maker. Suppose the stage payoff equals 1 during stages \( n^n \leq t < (n+1)^{n+1} \) for odd  \( n \), and equals 0 for even \( n \). The resulting payoff stream alternates between long runs of 0s and 1s, and its Cesàro averages oscillate between values close to 0 and 1 without converging. Hence, no robust maxmin strategy or robust equilibrium exists. This example illustrates that a stochastic game may have weak robust equilibria but not robust equilibria.

\subsection{Uniform equilibrium payoffs and the uniform value}

We recall the definition of the much-studied \emph{uniform equilibrium payoff}. Let \( \gamma_n:=\gamma_{u_n} \) denote the average payoff in the \( n \)-stage game:
\[
\gamma_n(\sigma)[z] := E^z_\sigma \left[ \frac{1}{n} \sum_{t=1}^n g_t \right] = E^z_\sigma u_n(g),
\]
where \( u_n \) denotes the \( n \)-stage Cesàro valuation.

\begin{defi}
A vector of payoffs \( v = (v^i[z])_{i \in I,\, z \in S} \) is called a \emph{uniform equilibrium payoff} of the stochastic game \( \Gamma = \langle I, S, A, g, q \rangle \) if, for every \( \varepsilon > 0 \), there exists a strategy profile \( \sigma_\varepsilon \) and an integer \( n_\varepsilon \) such that
\[
\gamma^i_n(\sigma_\varepsilon)[z] \geq v^i[z] - \varepsilon \quad \text{and} \quad \gamma^i_n(\sigma_\varepsilon^{-i},\tau^i)[z] \leq v^i[z] + \varepsilon
\]
for all \( n \geq n_\varepsilon \), \( i \in I \), \( z \in S \), and \( \tau^i \in \Sigma^i \).
\end{defi}

\begin{defi}\label{def:uniform-approximate-equilibrium-payoff} The stochastic game \(\Gamma\) has uniform approximate equilibrium
payoffs if, for every \(\varepsilon>0\), there exist a bounded payoff vector \(v_\varepsilon=(v^i_\varepsilon[z])_{i\in I,z\in S}\), a strategy profile \(\sigma_\varepsilon\), and an integer \(n_\varepsilon\) such that, for every \(n\ge n_\varepsilon\), every player \(i\), every initial state \(z\), and every unilateral deviation \(\tau^i\in\Sigma^i\),
\[
\gamma^i_n(\sigma_\varepsilon)[z]\ge v^i_\varepsilon[z]-\varepsilon
\]
and
\[
\gamma^i_n(\sigma_\varepsilon^{-i},\tau^i)[z]\le
v^i_\varepsilon[z]+\varepsilon .
\]
\end{defi}

Equivalently, since \(\varepsilon\) is arbitrary, \(\Gamma\) has uniform
approximate equilibrium payoffs if and only if, for every \(\varepsilon>0\), there exist a strategy profile \(\sigma_\varepsilon\) and an integer \(n_\varepsilon\) such that
\[
\gamma^i_n(\sigma_\varepsilon)[z]
\ge
\gamma^i_m(\sigma_\varepsilon^{-i},\tau^i)[z]-\varepsilon
\]
for all \(n,m\ge n_\varepsilon\), \(i\in I\), \(z\in S\), and
\(\tau^i\in\Sigma^i\).

\smallskip

When \(S\) is finite, this approximate notion is equivalent to the usual
existence of a uniform equilibrium payoff. Indeed, the usual notion immediately implies the approximate one; conversely, applying the approximate notion with \(\varepsilon_k\downarrow0\), every accumulation point of the bounded sequence of payoff vectors \(v_{\varepsilon_k}\in\mathbb R^{I\times S}\) is a uniform equilibrium payoff.

\smallskip
\noindent
\textit{Remark.} A uniform equilibrium payoff need not exist, and when it does, it may not be unique. By contrast, in two-player zero-sum stochastic games, the uniform equilibrium payoff---if it exists---is unique. In two-player zero-sum finite stochastic games, the uniform equilibrium payoff always exists, whereas in finite nonzero-sum stochastic games, the existence of a uniform equilibrium payoff remains an open problem.

In a two-player zero-sum stochastic game, the \emph{uniform value} is the uniform equilibrium payoff of Player~1, if it exists.

\subsection{Relation to the uniform value}In two-player zero-sum stochastic games, whenever the usual uniform value exists, it coincides with the uniform maxmin value of Player~1.

If \(v^i\) is the uniform maxmin value of player~\(i\), then
\[
\lim_{n\to\infty} v^i_n[z]=v^i[z]
\quad\text{for every }z\in S.
\]
In particular, the uniform maxmin value, if it exists, is unique.

\subsection{Sunspot, correlated, and extensive-form correlated strategy profiles}
An extensive-form correlated strategy profile is a mediated strategy profile in which a mediation device sends private signals to the players over time. At each stage, a player's action may depend not only on the public history of play but also on the private signals received by that player up to that stage.

Two important cases are sunspot and correlated strategy profiles. In a sunspot strategy profile, before play begins, a public signal is drawn and observed by all players. In a correlated strategy profile, before play begins, private signals are drawn and each player observes only his own signal. Thus, sunspot and correlated strategy profiles may be viewed as special classes of extensive-form correlated strategy profiles.

The classes considered in this paper are closed under mixing: a convex combination is implemented by first drawing an auxiliary randomization before play begins and then using the corresponding strategy profile. This convexity property will be used below to construct uniformly continuous selections by means of partitions of unity.

Let \(\mathcal M\) be a class of mediated strategy profiles. An
\(\mathcal M\)-strategy profile is a mediated strategy profile belonging to \(\mathcal M\). An \(\mathcal M\) \(\varepsilon\)-equilibrium is an
\(\mathcal M\)-strategy profile satisfying the corresponding
\(\varepsilon\)-incentive constraints against admissible unilateral deviations. The corresponding robustness and uniform payoff notions are defined by replacing ordinary strategy profiles by \(\mathcal M\)-strategy profiles, and ordinary unilateral deviations by admissible unilateral deviations in the corresponding mediated environment.

Thus, an \(\mathcal M\)-strategy profile \(\sigma\) is a robust \(\mathcal M\) \(\varepsilon\)-equilibrium for \(U\subseteq W^I\) if, for every \(i\in I\), \(u,w\in U\), initial state \(z\in S\), and admissible unilateral deviation \(\tau^i\) by player~\(i\),
\[
\gamma^i_u(\sigma^{-i},\tau^i)[z]\le \gamma^i_w(\sigma)[z]+\varepsilon.
\]

A vector \(v=(v^i[z])_{i\in I,z\in S}\) is a uniform \(\mathcal M\) equilibrium payoff if, for every \(\varepsilon>0\), there exist an \(\mathcal M\)-strategy profile \(\sigma_\varepsilon\) and an integer \(n_\varepsilon\) such that, for every \(n\ge n_\varepsilon\), every player \(i\), every initial state \(z\), and every admissible unilateral deviation \(\tau^i\),
\[
\gamma^i_n(\sigma_\varepsilon)[z]\ge v^i[z]-\varepsilon
\quad\text{and}\quad
\gamma^i_n(\sigma_\varepsilon^{-i},\tau^i)[z]\le v^i[z]+\varepsilon.
\]

Similarly, \(\Gamma\) has \emph{uniform approximate \(\mathcal M\) equilibrium payoffs} if, for every \(\varepsilon>0\), there exist a bounded payoff vector \(v_\varepsilon=(v^i_\varepsilon[z])_{i\in I,z\in S}\), an \(\mathcal M\)-strategy profile \(\sigma_\varepsilon\), and an integer \(n_\varepsilon\) such that, for every \(n\ge n_\varepsilon\), every player \(i\), every initial state \(z\), and every admissible unilateral deviation \(\tau^i\),
\[
\gamma^i_n(\sigma_\varepsilon)[z]\ge v^i_\varepsilon[z]-\varepsilon
\]
and
\[
\gamma^i_n(\sigma_\varepsilon^{-i},\tau^i)[z]\le
v^i_\varepsilon[z]+\varepsilon .
\]

\section{The main results}\label{sec:results}
\paragraph{\textbf{Overview.}}
Although our results span several settings—zero-sum and non-zero-sum, finite and infinite, two players and many—they follow from one principle, refined along a single axis. For each solution concept (the maxmin value and optimal strategies, Nash equilibrium, and the sunspot, correlated, and extensive-form correlated equilibria), robustness to small imprecision in the impatient valuations holds if and only if the appropriate uniform long-run condition holds. In the adversary case this condition is the existence of the uniform maxmin
value; in the Nash and mediated cases it is the existence of uniform approximate equilibrium payoffs, which, for finite state spaces, is equivalent to the corresponding usual exact uniform equilibrium-payoff notion. The form of this robustness is governed by convexity:
for the convex-valued concepts---maxmin strategies in the adversary case and mediated equilibrium strategy profiles---one can select uniformly continuously, over the entire valuation space, strategies satisfying the corresponding robust \(\varepsilon\)-maxmin or \(\varepsilon\)-equilibrium conditions, whereas for Nash equilibrium, whose equilibrium set is not convex, only a local, finite-cover selection is established.
The remaining distinctions---the cardinality of the state and action sets, and the number of players---do not alter this picture;
they determine only when the relevant uniform long-run condition is known to hold, through classical theorems (\cite{merney81} for finite two-person zero-sum games; \cite{ney03minmax} for finite stochastic games; \cite{vie00a,vie00b} for finite two-player stochastic games; \cite{solan99} for three-player absorbing games; \cite{solvie02} for extensive-form correlated equilibria). Tables~\ref{tab:robustness_comparison} and~\ref{tab:mediated_summary} summarize which guarantee applies in each case.

Our main result in the adversary case establishes the existence of a uniformly continuous selection of robust \(\varepsilon\)-maxmin strategies.\footnote{All continuity statements for strategy selections below use the metric convention for strategy spaces introduced in Section~3.2.}
Specifically, we show that in any stochastic game \(\Gamma\) in which the uniform maxmin value of player~\(i\) exists (e.g., any finite stochastic game), for every \(\varepsilon>0\) there exist \(\delta>0\) and a uniformly continuous mapping
\[
w\mapsto \sigma^i_w
\]
from \(W\) to the set of strategies of player~\(i\), such that \(\sigma^i_w\) is a robust \(\varepsilon\)-maxmin strategy of player~\(i\) for \(B(w,\delta)\cap W\).
As a consequence, the maxmin value mapping
\[
w\mapsto v^i(w)
\]
is uniformly continuous.

A key technical step in the proof of this result---and of the corresponding nonzero-sum results stated below---is the compactness of the completion of the metric space of impatient valuations. This compactness property, which is of independent interest, plays a central role in our analysis and may also be useful in the analysis of robust optimization in other dynamic models involving impatient valuations.

The uniform continuity of the maxmin value as a function of the impatient valuation in finite stochastic games follows directly from the existence of robust maxmin strategies or from the strategy selection conclusion itself.
Notably, this implication does not rely on auxiliary tools such as compactness of the completion of the metric space of impatient valuations or the existence of the uniform maxmin but follows directly from either the existence of robust maxmin strategies or a uniformly continuous selection of $\varepsilon$-maxmin strategies.

An alternative proof of the uniform continuity of the maxmin value in finite stochastic games can be obtained by combining two ingredients: the compactness of the completion of the metric space of impatient valuations, which is established in this paper, and the convergence of $v^i(w)$ as $w_1 \to 0$, which was established by Ziliotto \citeyearpar{zil2016,zil2018}.
While Ziliotto’s work does not directly address uniform continuity, it provides a key asymptotic result. Combined with the compactness property established here, this result yields the desired continuity conclusion.

Our result on the uniform continuity of the maxmin value is of independent interest, and its connection to prior asymptotic results further illustrates its scope and strength. It strictly strengthens  Ziliotto's \citeyearpar{zil2016,zil2018} convergence result, which shows that \( v^i(w) \) converges as \( w_1 \to 0 \). Ziliotto’s result, in turn, strictly strengthens the classical theorem of Bewley and Kohlberg \citeyearpar{bewkoh76}, which establishes that the limits of the \( n \)-stage and the \( \lambda \)-discounted values, as \( n \to \infty \) and \( \lambda \to 0^+ \), respectively, exist and coincide. Thus, our uniform continuity result unifies and extends the asymptotic convergence properties identified in both earlier contributions.

Importantly, our result extends beyond finite games: it establishes the uniform continuity of the maxmin value for any stochastic game that admits a uniform maxmin value—not just for finite games—thus extending the scope of both Ziliotto's and Bewley and Kohlberg's results.
\begin{remark}\label{rem:L1_payoff_continuity}
If the stage payoffs are uniformly bounded by \( C \), then for every strategy profile \( \sigma \) and any two impatient valuations \( w, w' \in W \), we have
\[
|\gamma^i_w(\sigma) - \gamma^i_{w'}(\sigma)| \leq C \sum_{t=1}^\infty |w_t - w'_t| = C \|w - w'\|_1,
\]
where \( \| \cdot \|_1 \) denotes the \( \ell_1 \) norm. In particular, the mapping \( w \mapsto v^i(w) \) is Lipschitz continuous (and therefore continuous) when \( W \) is equipped with the \( \| \cdot \|_1 \) metric.

Since each sequence \( w \in W \) consists of nonnegative terms summing to 1—and hence converging to 0—the \( \ell_1 \) and \( \ell_\infty \) norms induce the same topology on \( W \). Therefore, the mapping \( w \mapsto v^i(w) \) is also continuous under the topology induced by the \( \ell_\infty \) norm.

However, since \( W \) is not compact in the \( \ell_\infty \) norm, continuity alone does not imply uniform continuity—and it is this limitation that motivates our extension to the completion of \( W \).
\end{remark}

The assumptions underlying our main result in the adversary case vary across classes of stochastic games. In finite stochastic games, the existence of the uniform maxmin value is guaranteed by the result of \cite{ney03minmax}, building on \cite{merney81}. In general stochastic games, however, this existence must be assumed explicitly, as it does not hold universally.

In the nonzero-sum case, our main result establishes a weaker form of regularity for \(\varepsilon\)-equilibrium strategies. We show, more generally, that in any stochastic game \(\Gamma\) that has uniform approximate equilibrium payoffs---and hence, in particular, in any game that has a uniform equilibrium payoff, such as any finite two-player stochastic game, any finite irreducible stochastic game, any finite stochastic game with perfect information, or any three-player absorbing game---there exist, for every \(\varepsilon>0\), a finite cover \((U_j)_{j\in J}\) of the space \(W^I\) of impatient valuation profiles, a corresponding family of strategy profiles \((\sigma_j)_{j\in J}\),
and a constant \(\delta>0\), such that, for each \(j\in J\), the strategy profile \(\sigma_j\) is a robust \(\varepsilon\)-equilibrium for the \(\delta\)-neighborhood of \(U_j\).\footnote{By definition, a robust \(\varepsilon\)-equilibrium for a set of profiles of impatient valuations implies that \(\sigma_j\) is an \(\varepsilon\)-equilibrium for \(U_j^\delta\) and that the \(\sigma_j\)-payoff varies by at most \(\varepsilon\) over \(U_j^\delta\).}

The weaker regularity in the nonzero-sum case contrasts with the uniform continuity obtained in the adversary case. This distinction stems from the geometry of the solution sets: the set of $\varepsilon$-maxmin strategies is convex, enabling the use of partition-of-unity arguments over the compact completion of $W$. By contrast, the set of $\varepsilon$-equilibrium profiles is typically nonconvex, hence, a mixture of $\varepsilon$-equilibrium profiles need not be an $\varepsilon$-equilibrium, which precludes the use of partition-of-unity arguments to construct $\varepsilon$-equilibrium profiles as mixtures of $\varepsilon$-equilibrium profiles.

The situation changes when the equilibrium concept is convex-valued: stronger results become available, as in the case of sunspot $\varepsilon$-equilibria and extensive-form correlated (EFC) $\varepsilon$-equilibria. Convexity, together with the compactness of the completion of the metric space of impatient valuations, enables the use of partition-of-unity techniques to obtain uniform continuity results across the entire domain.

For sunspot, correlated, and extensive-form correlated equilibria, these
regularity conclusions are part of full characterizations: the existence of uniform approximate mediated equilibrium payoffs is equivalent to the existence of robust mediated equilibria, and also to the existence of uniformly continuous robust approximate selections. In finite state spaces, this is equivalent to the usual existence of a uniform mediated equilibrium payoff.

Although this paper focuses on discrete-time stochastic games, the underlying techniques extend to other multistage environments, including continuous-time stochastic games. In particular, the robustness approach developed above extends to finite non-zero-sum continuous-time settings: for every finite continuous-time stochastic game $\Gamma$ and every $\varepsilon>0$, there exists a uniformly continuous mapping $w\mapsto \sigma_w^\varepsilon$ from \(W^I\), the space of profiles of impatient valuations, to strategy profiles, such that $\sigma_w^\varepsilon$ is an $\varepsilon$-equilibrium of $\Gamma_w$.\footnote{The proof relies on the existence of a uniform equilibrium payoff in finite continuous-time stochastic games~\cite{neyman2017}, and the mixing of $\varepsilon$-equilibrium uses techniques of jointly controlled lotteries, which are ``costless'' in continuous-time games.}

\subsection*{Formal Statements of the Main Results}

Before presenting the general result for the adversary case, we begin with two corollaries that apply to finite stochastic games and establish the uniform continuity of the maxmin value and the existence, for every \(\varepsilon>0\), of a uniformly continuous selection of robust \(\varepsilon\)-maxmin strategies.

The sup-norm on the set  of profiles of impatient valuations is defined as follows. For \(w=(w^i)_{i\in I}, u=(u^i)_{i\in I}\in W^I\), we define \[\|w-u\|:=\max_{i\in I}\|w^i-u^i\|_\infty.\]

\begin{thm}[{\bf Uniform continuity of the maxmin value function}]\label{thm:uniformcontinuityofv}
Let $\Gamma$ be a finite stochastic game. Then the mapping \(w \mapsto v(w)\), from \(W^I\) to \(\mathbb{R}^{I\times S}\), is uniformly continuous when \(W^I\) is equipped with the sup-norm.
\end{thm}

\begin{thm}[{\bf A uniformly continuous selection of robust \(\varepsilon\)-maxmin strategies}]
\label{thm:uniformcontinuousoptimalstrategies}
Let $\Gamma$ be a finite stochastic game.
Then, for every \(\varepsilon>0\) and every player~\(i\), there exist \(\delta>0\) and a uniformly continuous mapping
\[
w^i\mapsto \sigma^i_{w^i,\varepsilon}
\]
from \(W\) to \(\Sigma^i\), such that, for every \(w^i\in W\), the strategy
\(\sigma^i_{w^i,\varepsilon}\) is a robust \(\varepsilon\)-maxmin strategy of player~\(i\) for
\[
B(w^i,\delta)\cap W .
\]
\end{thm}

\begin{remark}
Since the space of mixed strategies (endowed with the total variation distance) is a complete metric space, and the mapping $w\mapsto \sigma^i_{w^i,\varepsilon}$ is uniformly continuous, the limit
$\sigma^i_\varepsilon:=\lim_{w^i_1 \to 0}\sigma^{i}_{w^i,\varepsilon}$ exists. Using the uniform continuity of
$w^i \mapsto v^i(w^i)$, it follows that
there exists \(\delta>0\) such that \(\sigma^i_\varepsilon\) is a robust \(2\varepsilon\)-maxmin strategy of player~\(i\) for the set
\[
\{u^i\in W:u^i_1<\delta\}.
\]
The last conclusion recovers the classical uniform-value/maxmin conclusion in robust terminology. In the two-player zero-sum case, the existence of a uniform value implies that, for every \(\varepsilon>0\), there exist a strategy \(\sigma_\varepsilon\) and an integer \(k_\varepsilon\) such that \(\sigma_\varepsilon\) is a robust \(\varepsilon\)-maxmin strategy for the set \[ \{u_n:n\ge k_\varepsilon\}. \] Equivalently, using Lemma~\ref{lem:cesaro_patient}, there exist a strategy \(\sigma_\varepsilon\) and \(\delta_\varepsilon>0\) such that \(\sigma_\varepsilon\) is a robust \(\varepsilon\)-maxmin strategy for \[ \{w\in W:w_1\le \delta_\varepsilon\}. \]
Thus Theorem~\ref{thm:uniformcontinuousoptimalstrategies} may be viewed as a uniformly continuous selection strengthening of the classical uniform value and uniform maxmin results of \citet{merney81} and \citet{ney03minmax}.
\end{remark}

\smallskip
We now state the general result from which the preceding two theorems follow. It characterizes the existence of robust \(\varepsilon\)-maxmin strategies that vary uniformly continuously with the impatient valuation.
The result applies to general (not necessarily finite) stochastic games, under the assumption that the uniform maxmin value exists.

\begin{thm}[{\bf Characterization of Stochastic Games Having Robust Maxmin Strategies}]\label{thm:robustoptimalstrategy}
Let \(\Gamma\) be a stochastic game. Then the following conditions are equivalent for a fixed player~\(i\):
\begin{conditions}
\item[\((1.a)\)] \(\Gamma\) has a uniform maxmin value for player~\(i\).
\item[\((1.b)\)] \(\Gamma\) has robust maxmin strategies for player~\(i\).
\item[\((1.b^\ast)\)] For every \(\varepsilon>0\), there exists \(\delta>0\) such that, for every \(U\subseteq W\) with \(\operatorname{diam}(U)<\delta\), player~\(i\) has a robust \(\varepsilon\)-maxmin strategy for \(U\). This is the definition
of \((1.b)\).
\item[\((1.c)\)] For every \(\varepsilon > 0\), there exist \(\delta > 0\) and a finite family \((\sigma_j, U_j)_{j \in J}\), where \(\sigma_j\) is a strategy of player~\(i\) and \(U_j \subseteq W\), such that:
    \begin{itemize}
    \item[] \(\bigcup_{j \in J} U_j = W\), and
    \item[] \(\sigma_j\) is a robust \(\varepsilon\)-maxmin strategy of player~\(i\) for \[ U_j^\delta := \{w' \in W : d(w', U_j) < \delta\}. \]
        \end{itemize}
\item[\((1.d)\)] For every \(\varepsilon > 0\), there exist \(\delta > 0\) and a uniformly continuous mapping \(w \mapsto \sigma_w\) from \(W\) to the set of strategies of player~\(i\), such that \(\sigma_w\) is a robust \(\varepsilon\)-maxmin strategy of player~\(i\) for \(B(w,\delta)\cap W\).
    \end{conditions}
\end{thm}

\begin{remark}
In finite two-player zero-sum stochastic games—or even finite MDPs—the mapping from impatient valuations to optimal strategies need not be continuous.
For example, consider an MDP $\Gamma$ with one nonabsorbing state and two actions, and two absorbing states with absorbing payoffs of 0 and 1. In the nonabsorbing state, one action yields a stage payoff of 1 and transitions to the absorbing state with payoff 0, while the other action yields a stage payoff of 0 and transitions to the absorbing state with payoff 1. If $w_1 < 1/2$, the only optimal strategy is to play the action that yields a stage payoff of 0; if $w_1 > 1/2$, the only optimal strategy is to play the action that yields a stage payoff of 1. Thus, the optimal strategy switches discontinuously at $w_1 = 1/2$.

The discontinuity in this example is not an artifact of the particular metric placed on strategies: the unique optimal strategies on the two sides of \(w_1=1/2\) are distinct, so the optimal-strategy selection is discontinuous under any metric that separates distinct strategies. By contrast, Theorem~\ref{thm:uniformcontinuousoptimalstrategies} guarantees the existence of a uniformly continuous mapping from impatient valuations to approximate optimal strategies in the total variation metric on mixed strategies. This total variation continuity implies continuity for the standard weaker behavioral topologies and metrics on strategies. Moreover, the uniformly continuous selection admits a limiting strategy
\[
\sigma_\varepsilon:=\lim_{w_1\to0}\sigma_{w,\varepsilon},
\]
and for every \(\eta>\varepsilon\), there exists \(\delta_\eta>0\) such that
\(\sigma_\varepsilon\) is \(\eta\)-optimal in \(\Gamma_w\) for every
\(w\in W\) with \(w_1<\delta_\eta\).
\end{remark}

\begin{remark}\label{rem:condition-1-a-star}
We introduce the following property, which is implied by \((1.a)\).

\((1.a^\ast)\) \(\Gamma\) has uniform approximate maxmin payoffs for player~\(i\): for every \(\varepsilon>0\), there exist a bounded vector
\(v^i_\varepsilon=(v^i_\varepsilon[z])_{z\in S}\), an integer \(n_\varepsilon\), and a strategy \(\sigma^i_\varepsilon\in\Sigma^i\) such that, for every \(n\ge n_\varepsilon\),
\[
\gamma^i_n(\sigma^i_\varepsilon,\tau)[z]
\ge v^i_\varepsilon[z]-\varepsilon
\quad\text{for all }\tau\in\Sigma^{-i},\ z\in S,
\]
and, for every \(\rho^i\in\Sigma^i\) and every \(w\in W\) with
\(w_1<1/n_\varepsilon\), there exists \(\tau=\tau(\rho^i,w)\in\Sigma^{-i}\) such that
\[
\gamma^i_w(\rho^i,\tau)[z]
\le v^i_\varepsilon[z]+\varepsilon
\quad\text{for all }z\in S.
\]

For maxmin payoffs, \((1.a)\) and \((1.a^\ast)\) are equivalent, even when \(S\) is infinite. Indeed, if \(v^i_\varepsilon\) and \(v^i_\eta\) satisfy the two approximate conditions with tolerances \(\varepsilon\) and \(\eta\), then applying the lower condition for one vector and the upper condition for the other to sufficiently long Cesaro valuations gives
\[
\|v^i_\varepsilon-v^i_\eta\|_\infty\le \varepsilon+\eta .
\]
Thus, whenever \(\varepsilon_j\downarrow 0\), the sequence
\((v^i_{\varepsilon_j})_j\) is Cauchy in \(\ell^\infty(S)\), and hence converges in sup norm. The estimate above also shows that the limit is independent of the chosen sequence \(\varepsilon_j\downarrow 0\). Passing to this limit in the two approximate conditions shows that the limit is the uniform maxmin value.
\end{remark}

We next consider the mediated equilibrium concepts defined in Section~3.8: sunspot, correlated, and extensive-form correlated equilibria. For these classes, the closed-under-mixing structure allows the local robustness conclusions to be strengthened to uniformly continuous robust selections, in parallel with the maxmin case.

For finite stochastic games, the existence of a uniform extensive-form correlated equilibrium payoff is known; hence the characterization below implies the existence of uniformly continuous robust extensive-form correlated \(\varepsilon\)-equilibrium selections in every finite stochastic game.

\begin{thm}[{\bf Characterization of robust mediated equilibria}]
\label{thm:characterisation_mediated_equilibria}
Let \(\Gamma\) be an \(n\)-player stochastic game, and let \(\mathcal M\) be one of the following classes of mediated strategy profiles: sunspot, correlated, or extensive-form correlated strategy profiles. Then the
following conditions are equivalent:
\begin{conditions}
\item[\((M.a^\ast)\)] \(\Gamma\) has uniform approximate \(\mathcal M\) equilibrium payoffs.

\item[\((M.b)\)] \(\Gamma\) has robust \(\mathcal M\) equilibria.

\item[\((M.b^\ast)\)] For every \(\varepsilon>0\), there exists \(\delta>0\) such that, for every \(U\subseteq W^I\) with \(\operatorname{diam}(U)<\delta\), there exists a robust \(\mathcal M\) \(\varepsilon\)-equilibrium for \(U\). This is the definition of \((M.b)\).

\item[\((M.c)\)] For every \(\varepsilon>0\), there exist \(\delta>0\) and a finite family \((\sigma_j,U_j)_{j\in J}\), where each \(\sigma_j\) is an \(\mathcal M\)-strategy profile and \(U_j\subseteq W^I\), such that
\[
\bigcup_{j\in J}U_j=W^I,
\]
and \(\sigma_j\) is a robust \(\mathcal M\) \(\varepsilon\)-equilibrium for
\[
U_j^\delta:=\{w'\in W^I:d(w',U_j)<\delta\}.
\]

\item[\((M.d)\)] For every \(\varepsilon>0\), there exist \(\delta>0\) and a uniformly continuous mapping \(w\mapsto\sigma_w\) from \(W^I\) to the set of \(\mathcal M\)-strategy profiles, such that, for every \(w\in W^I\), the profile \(\sigma_w\) is a robust \(\mathcal M\) \(\varepsilon\)-equilibrium for
\[
B(w,\delta)\cap W^I.
\]
\end{conditions}
\end{thm}

\begin{remark}
If \(S\) is finite, then condition \((M.a^\ast)\) is equivalent to the existence of a uniform \(\mathcal M\) equilibrium payoff in the usual sense. The implication from the usual notion to \((M.a^\ast)\) is immediate. Conversely, applying \((M.a^\ast)\) with \(\varepsilon_k\downarrow0\), every accumulation point of the bounded sequence of  payoff vectors \(v_{\varepsilon_k}\in\mathbb R^{I\times S}\) is a uniform \(\mathcal M\) equilibrium payoff.

For \(\mathcal M\) equal to the class of extensive-form correlated strategy profiles, condition \((M.a^\ast)\) is satisfied in every finite stochastic game. For \(\mathcal M\) equal to the class of sunspot strategy profiles, the existence of a uniform equilibrium payoff implies condition \((M.a^\ast)\).

\end{remark}

\begin{remark}
The stronger assumption that the stochastic game admits a uniform equilibrium payoff is satisfied in several important classes of stochastic games, including:
\begin{itemize}
\item all finite two-player stochastic games,
\item finite irreducible stochastic games,
\item finite stochastic games with perfect information, and
\item three-player absorbing games.
\end{itemize}
\end{remark}

\begin{remark}
The assumption that the game admits a uniform sunspot equilibrium payoff is also satisfied in several special classes of absorbing games, including those with a single nonabsorbing action profile~\cite{solan2020quitting}, those in which all nonabsorbing payoffs equal 0 and all absorbing payoffs are strictly positive, those with at most three nonabsorbing action profiles~\cite{munk2020sunspot}, and those in which the nonabsorbing action profiles form a multi-dimensional ``submatrix''~\cite{solan2021sunspot}.
\end{remark}

We now present the general result for the nonzero-sum case, which provides a characterization of robust approximate equilibria in arbitrary stochastic games that have uniform approximate equilibrium payoffs, without requiring convexity of the equilibrium set. In finite state spaces, this condition is equivalent to the existence of a uniform equilibrium payoff.

\begin{thm}[{\bf Characterization of Stochastic Games Having Robust Equilibria}]\label{thm:robustequilibrium}
The following conditions on an $n$-player stochastic game $\Gamma$ are equivalent:
\begin{conditions}
\item[\((2.a^\ast)\)] $\Gamma$ has uniform approximate equilibrium payoffs.
\item[\((2.b)\)] $\Gamma$ has robust equilibria.
\item[\((2.b*)\)] For every \(\varepsilon>0\), there exists \(\delta>0\) such that, for every \(U\subseteq W^I\) with \(\operatorname{diam}(U)<\delta\), there exists a robust \(\varepsilon\)-equilibrium for \(U\). This is the definition of \((2.b)\).
\item[\((2.c)\)] For every $\varepsilon > 0$, there exist $\delta > 0$ and a finite family $(\sigma_j, U_j)_{j \in J}$, where each $\sigma_j$ is a strategy profile and $U_j \subseteq W^I$, such that:
\begin{itemize}
\item[] $\bigcup_{j \in J} U_j = W^I$, and
\item[] \(\sigma_j\) is a robust \(\varepsilon\)-equilibrium for
\[
U_j^\delta:=\{w'\in W^I:d(w',U_j)<\delta\}.
\]
\end{itemize}
\end{conditions}
\end{thm}

\begin{remark}
If \(S\) is finite, then condition \((2.a^\ast)\) is equivalent to the existence of a uniform equilibrium payoff in the usual sense. The implication from the usual notion to \((2.a^\ast)\) is immediate. Conversely, applying \((2.a^\ast)\) with \(\varepsilon_k\downarrow0\), every accumulation point of the bounded sequence of payoff vectors \(v_{\varepsilon_k}\in\mathbb R^{I\times S}\) is a uniform equilibrium payoff.
\end{remark}

\begin{remark}A closely related finite-valued formulation of \((2.c)\), equivalent after a harmless change of \(\varepsilon\), is the following:

\begin{itemize}
\item[\((2.c^\ast)\):]
For every \(\varepsilon>0\), there exist \(\delta>0\) and a
finite-valued mapping
\[
W^I\ni w\mapsto(\sigma_w,v(w))\in\Sigma\times\mathbb R^{I\times S}
\]
such that, for every \(w'\in B(w,\delta)\cap W^I\), the profile \(\sigma_w\) is an \(\varepsilon\)-equilibrium in \(\Gamma_{w'}\), and
\[
\|\gamma_{w'}(\sigma_w)-v(w)\|<\varepsilon .
\]
\end{itemize}

Indeed, if (2.c) holds, we may assume, without loss of generality, that the sets \( U_j \) are disjoint. Define \( \sigma_w = \sigma_j \) if \( w \in U_j \), and define \( v(w) \) to be any point in \( \{ \gamma_u(\sigma_j) : u \in U_j \} \).

Conversely, if \((2.c^\ast)\) holds, define the sets \(U_j\) as the nonempty pre-images of the mapping \(W^I\ni w\mapsto(\sigma_w,v(w))\).
\end{remark}

\begin{remark}
Unlike Theorem~3 in the adversary case and Theorem~4 in the mediated case, Theorem~5 does not assert a global uniformly continuous selection of robust approximate Nash equilibria, nor a global continuous selection of approximate equilibrium payoffs. The reason is that the proof of the corresponding selection results uses convex mixtures of local solutions. This argument is unavailable for Nash equilibrium, since Nash \(\varepsilon\)-equilibrium profiles are generally nonconvex.

\end{remark}

\section{Summary of Results}\label{sec:summary}

The two tables below summarize the main robustness and continuity results of the paper, for the adversary, Nash, and mediated-equilibrium cases. Together they trace a single chain: the appropriate uniform long-run condition is equivalent to robustness. In the adversary case this condition is the existence of the uniform maxmin value;
in the Nash and mediated cases it is the existence of uniform approximate equilibrium payoffs, which, for finite state spaces, is equivalent to the corresponding usual exact uniform equilibrium-payoff notion.

The resulting robustness strengthens from a local, finite-cover selection for Nash equilibrium to a global, uniformly continuous selection for the convex-valued concepts: maxmin strategies in the adversary case, and sunspot, correlated, and extensive-form correlated equilibria in
the mediated case. We split the summary into two tables only for legibility. Recall that \(W^I\) denotes the space of impatient valuation profiles, and let \(\Gamma_w\) be the stochastic game defined by the impatient valuation profile \(w\).

We recall that a strategy \(\sigma^i\) is a robust \(\varepsilon\)-maxmin strategy of player~\(i\) for a set \(U\subseteq W\) if it satisfies the robust maxmin inequality of Definition~\ref{def:robust_maxmin_strategy_for_U} for all valuations in \(U\). Similarly, a profile \(\sigma\) is a robust \(\varepsilon\)-equilibrium for \(U\subseteq W^I\) if it satisfies the robust equilibrium inequality of Definition~\ref{def:robust_equilibrium_for_U} for all valuation profiles in \(U\). These conditions imply, in particular, the corresponding pointwise \(\varepsilon\)-optimality or \(\varepsilon\)-equilibrium conditions and payoff stability across \(U\).
For \(\delta>0\), \(U^\delta\) denotes the \(\delta\)-neighborhood of \(U\), with respect to the sup-norm on \(W\) or \(W^I\), as appropriate.

In the adversary column, the statements are for a fixed player~\(i\). Thus \(w\in W\) denotes player~\(i\)'s impatient valuation, and \(v^i(w)\in\mathbb{R}^S\) denotes the maxmin payoff vector of player~\(i\) at valuation \(w\).

\begin{table}[H]
\centering
\renewcommand{\arraystretch}{1.3}
\begin{tabular}{|p{7.9cm}|p{7.9cm}|}
\hline
\textbf{Maxmin/Adversary Case} & \textbf{Non-zero-Sum Case} \\
\hline\hline
{Classical uniform property: \(\Gamma\) has a uniform maxmin value for player~\(i\).}
&%
{Classical uniform property: \(\Gamma\) has a uniform equilibrium payoff.} \\
\hline
\(\Longleftrightarrow\ \Gamma\) has uniform approximate maxmin payoffs for player~\(i\).
& \(\Longrightarrow\) \(\Gamma\) has uniform approximate equilibrium payoffs.
For finite \(S\), the converse implication also holds.\\
\hline
\(\Longleftrightarrow\ \Gamma\) has robust maxmin strategies for player~\(i\) (Thm.~\ref{thm:robustoptimalstrategy}).
& \(\Longleftrightarrow\  \Gamma\) has robust equilibria (Thm.~\ref{thm:robustequilibrium}).\\
\hline
\(\Longleftrightarrow\ \forall\varepsilon>0,\ \exists\delta>0,\ \exists\) a finite cover \((U_j)\) of \(W\) and strategies \((\sigma_j)\) of player~\(i\) such that \(\sigma_j\) is a robust \(\varepsilon\)-maxmin strategy for \(U_j^\delta\).
& \(\Longleftrightarrow\ \forall \varepsilon>0\), \(\exists \delta>0\),  a finite cover \((U_j)\) of \(W^I\), and strategy profiles \((\sigma_j)\) such that \(\sigma_j\) is a robust \(\varepsilon\)-equilibrium for \(U_j^{\delta}\). \\
\hline
\(\Longleftrightarrow\)  \(\forall \varepsilon>0\), \(\exists\) \(\delta>0\) and a uniformly continuous map
\(W\ni w\mapsto\sigma_w\in\Sigma^i\) such that \(\sigma_w\) is a robust
\(\varepsilon\)-maxmin strategy for \(B(w,\delta)\cap W\).
& No global continuous selection result is asserted in general; only the finite-cover local robustness statement is proved.
\\
\hline
\(\Longrightarrow\) the mapping \(W\ni w\mapsto v^i(w)\in\mathbb{R}^S\) is uniformly continuous.
& No general global continuous payoff selection is asserted; on each patch \(U_j^\delta\), the fixed robust \(\varepsilon\)-equilibrium \(\sigma_j\) has payoff oscillation at most \(\varepsilon\).
 \\
\hline\hline
\end{tabular}
\caption{Comparison of robustness and continuity results in the adversary and non-zero-sum cases}
\label{tab:robustness_comparison}
\end{table}

\medskip
The preceding table focuses on robustness and continuity properties of maxmin and equilibrium strategies. The next table summarizes analogous results for sunspot, correlated, and extensive-form correlated (EFC) equilibria. In contrast to the nonconvexity of the Nash equilibrium correspondence, these correlated equilibrium concepts are convex-valued, which enables stronger regularity properties and the application of partition-of-unity techniques.

The existence of a uniform equilibrium payoff implies the existence of a uniform sunspot equilibrium payoff, which implies the existence of a uniform correlated equilibrium payoff, which in turn implies the existence of a uniform extensive-form correlated equilibrium payoff. Thus the existence assumptions below become weaker as the mediation device becomes richer.

\begin{table}[H]
\centering
\renewcommand{\arraystretch}{1.3}
\begin{tabular}{|p{15.8cm}|}
\hline
\textbf{Mediated equilibria: sunspot, correlated, or extensive-form correlated} \\
For \(\mathcal M\) equal to one of these three classes\\
\hline\hline
Classical uniform property: \(\Gamma\) has a uniform \(\mathcal M\) equilibrium payoff. \\
\hline
\(\Longrightarrow\ \Gamma\) has uniform approximate \(\mathcal M\) equilibrium payoffs.  For finite \(S\), the converse implication also holds. \\
\hline
\(\Longleftrightarrow\ \Gamma\) has robust \(\mathcal M\) equilibria. \\
\hline
\(\Longleftrightarrow\ \forall \varepsilon > 0\), \(\exists\ \delta>0\) and a uniformly continuous map
\(W^I \ni w \mapsto \sigma_w\) such that, for every \(w\in W^I\), \(\sigma_w\) is a robust \(\mathcal M\)
\(\varepsilon\)-equilibrium for \(B(w,\delta)\cap W^I\).\\
\hline
\end{tabular}
\caption{Continuity and robustness results for mediated correlation concepts}
\label{tab:mediated_summary}
\end{table}

\newpage
\section{Proof of the Theorems}\label{sec:proofofthms}
The first part of this section treats the adversarial case, in which the maxmin payoff and strategies of player~\(i\) depend only on player~\(i\)’s own valuation. A central step is to show that the existence of a uniform maxmin value yields, for every \(\varepsilon>0\), a uniformly continuous mapping \(w\mapsto\sigma^i_w\) that assigns to each impatient valuation \(w\in W\) an \(\varepsilon\)-maxmin strategy \(\sigma^i_w\) in \(\Gamma_w\). The construction exploits the compactness of the completion \(\overline W\) of the (noncompact) metric space \(W\) of impatient valuations. The main local step is to establish that, for every \(w\in\overline W\) and \(\varepsilon>0\), there exist \(\delta>0\) and a strategy \(\sigma^i\) such that \(\sigma^i\) is a robust \(\varepsilon\)-maxmin strategy of player~\(i\) for all impatient valuations \(u\in W\) with \(d(u,w)<\delta\).

Compactness yields a finite open cover \((U_j)_{j\in J}\) of \(\overline W\) and corresponding strategies \((\sigma^i_j)_{j\in J}\) such that, for each \(j\in J\), the strategy \(\sigma^i_j\) is a robust \(\varepsilon\)-maxmin strategy for \(U_j\cap W\). Using the convexity of the set of \(\varepsilon\)-maxmin mixed strategies and a partition of unity \((f_j)_{j\in J}\) \emph{subordinate} to \((U_j)_{j\in J}\) (nonnegative continuous functions with \(\sum_{j\in J} f_j(w)=1\) for all \(w\in\overline W\) and \(\operatorname{supp} f_j\subseteq U_j\)), we define
\[
w\longmapsto \sigma^i_w := \sum_{j\in J} f_j(w)\,\sigma^i_j.
\]
\(\sigma^i_w\) is the mixed strategy playing \(\sigma^i_j\) with probability \(f_j(w)\), and the maxmin guarantee survives mixing because \(\inf_\tau \sum_j f_j(w)y(\sigma^i_j,\tau)\geq \sum_j \inf_\tau f_j(w)y(\sigma^i_j,\tau)\).
By multilinearity of expected payoffs in mixed strategies, \(\sigma^i_w\) is \(\varepsilon\)-maxmin for every \(w\in W\); and since the sum is finite and each \(f_j\) is continuous on the compact domain \(\overline W\), the mapping \(w\mapsto \sigma^i_w\) is uniformly continuous.

We begin by introducing notation that will be used throughout the proofs.
Let \((\overline{W}, d)\) be the metric space of all nonincreasing sequences \(w=(w_t)_{t\ge1}\) with \(w_t\geq 0\) and \(\sum_{t=1}^\infty w_t\le1\), endowed with the supremum metric
\[
d(w,u):=\|w-u\|_\infty=\sup_{t\ge1}|w_t-u_t|.
\]
The inequality \(0\leq t\,w_t \le \sum_{s=1}^t w_s \le 1\) implies \(0\leq w_t \le 1/t\) for all \(t\); hence \(w_t\to 0\) as \(t\to\infty\).

Note that \(W\subset \overline{W}\) and the restriction of \(d\) to \(W\) coincides with the metric on \(W\) defined in Section~\ref{sec:basicconcepts}.

The space \(\overline W\) is compact, and \(W\) is dense in
\(\overline W\). Indeed, if \(w\in\overline W\) and
\[
\beta(w):=1-\sum_{t=1}^{\infty}w_t,
\]
then, for every \(n\),
\[
w+\beta(w)u_n\in W,
\]
and
\[
\|w+\beta(w)u_n-w\|_\infty
=
\beta(w)\|u_n\|_\infty
=
\frac{\beta(w)}{n}
\longrightarrow 0 .
\]

Let \((\overline{W}^I, d)\) denote the product metric space of profiles \(w=(w^i)_{i\in I}\), where each \(w^i\in\overline{W}\), with
\[
d(w,u):=\max_{i\in I} d\bigl(w^i,u^i\bigr).
\]

Next, we introduce two localized notions of robustness that will be used in the sequel.

\begin{defi}
Fix \( w \in \overline{W} \). A stochastic game is said to \emph{have robust maxmin strategies for player~\( i \) at \( w \)} if, for every \( \varepsilon > 0 \), there exist \( \delta > 0 \) and a strategy \( \sigma_\varepsilon \) of player~\( i \) such that \( \sigma_\varepsilon \) is a robust \( \varepsilon \)-maxmin strategy of player~\( i \) for the set \( \{u\in W : d(u,w)<\delta\} \).
\end{defi}

\begin{defi}
Fix \( w \in \overline{W}^I \). A stochastic game is said to \emph{have robust equilibria at \( w \)} if, for every \( \varepsilon > 0 \), there exist \( \delta > 0 \) and a strategy profile \( \sigma_\varepsilon \) such that \( \sigma_\varepsilon \) is a robust \( \varepsilon \)-equilibrium for the set \( \{u\in W^I : d(u,w)<\delta\} \).
\end{defi}

\smallskip

\noindent
\textit{Remark.} These notions express \emph{local} forms of robustness, defined at fixed points \(w\in\overline{W}\) or \(w\in\overline{W}^I\). In what follows, we first give a sufficient condition under which the stochastic game \(\Gamma\) has robust maxmin strategies at each \(w\in\overline{W}\), where the associated \(\delta\) may depend on \(w\). We then use the compactness of \(\overline{W}\) to show that, under this condition, \(\Gamma\) has robust maxmin strategies.

\subsection{Proof of Theorem~\ref{thm:robustoptimalstrategy}}\label{sec:proofofthm1}
The following elementary summation-by-parts lemma converts uniform bounds on sufficiently long Cesaro averages into uniform bounds for nonincreasing weights whose first weight is sufficiently small. It will be used both here and in the equilibrium proofs below.

\begin{lem}[Cesaro bounds and patient subvaluations]\label{lem:cesaro_patient}
Let \(M\ge 0\), let \(a\in\mathbb R\) satisfy \(|a|\le M\), let \(N\) be a positive integer, and let \(\eta>0\). Then, for every tolerance \(\varepsilon>\eta\), there exists \(\theta>0\) such that the following holds.

If \((x_t)_{t\ge1}\) is a real sequence with \(|x_t|\le M\) and
\[\frac1n\sum_{t=1}^n x_t\ge a-\eta
\quad\text{for every } n\ge N,\]
then, for every nonincreasing sequence \((\lambda_t)_{t\ge1}\) with
\(\lambda_t\ge0\), \(\sum_{t=1}^\infty\lambda_t\le1\), and \(\lambda_1<\theta\),
\[
\sum_{t=1}^\infty \lambda_t x_t
\ge
\left(\sum_{t=1}^\infty \lambda_t\right)a-\varepsilon .
\]

The analogous upper bound also holds: if
\[
\frac1n\sum_{t=1}^n x_t\le a+\eta
\quad\text{for every } n\ge N,
\]
then, for every such sequence \((\lambda_t)_{t\ge1}\) with \(\lambda_1<\theta\),
\[
\sum_{t=1}^\infty \lambda_t x_t
\le
\left(\sum_{t=1}^\infty \lambda_t\right)a+\varepsilon .
\]
\end{lem}

\begin{proof}
We prove the lower bound. The upper bound follows, with the same \(\theta\), by applying the lower bound to the sequence \((-x_t)_{t\ge1}\) and the constant \(-a\).

Choose \(\theta>0\) so small that
\[
(N-1)\theta(2M+\eta)<\varepsilon-\eta .
\]
This choice depends only on \(M,N,\eta,\varepsilon\) and not on the sequence \((x_t)\).

Now let \((x_t)_{t\ge1}\) and \((\lambda_t)_{t\ge1}\) satisfy the hypotheses of the lower-bound statement, and write
\[
A_n:=\frac{1}{n}\sum_{t=1}^n x_t .
\]
Since \((\lambda_t)_{t\ge 1}\) is nonincreasing and \(\lambda_t\to 0\), summation by parts gives
\[
\sum_{t=1}^\infty \lambda_t x_t
=
\sum_{n=1}^\infty n(\lambda_n-\lambda_{n+1})A_n .
\]
The coefficients \(n(\lambda_n-\lambda_{n+1})\) are nonnegative and sum to
\(\sum_{t=1}^\infty \lambda_t\). Moreover,
\[
\sum_{n=1}^{N-1} n(\lambda_n-\lambda_{n+1})
\le
\sum_{n=1}^{N-1} \lambda_n
\le
(N-1)\lambda_1 .
\]
Thus, if \(\lambda_1<\theta\), the total weight assigned to the terms \(n<N\) is at most
\((N-1)\theta\). For \(n\ge N\), we have \(A_n\ge a-\eta\), while for all \(n\), \(|A_n|\le M\). Hence
\[
\sum_{t=1}^\infty \lambda_t x_t
\ge
\left(\sum_{t=1}^\infty \lambda_t\right)(a-\eta)
-
(N-1)\theta(2M+\eta).
\]
Since \(\sum_{t=1}^\infty \lambda_t\le 1\), it follows that
\[
\sum_{t=1}^\infty \lambda_t x_t
\ge
\left(\sum_{t=1}^\infty \lambda_t\right)a
-\eta
-(N-1)\theta(2M+\eta)
\ge
\left(\sum_{t=1}^\infty \lambda_t\right)a-\varepsilon,
\]
where the last inequality follows from the choice of \(\theta\).
\end{proof}

We proceed with the following intermediate result, which is of independent interest. It shows that the existence of a uniform maxmin value implies local robustness.

\begin{lem}\label{lem:uniformmaxminatwinoverlineW}
Let \(\Gamma\) be a stochastic game that has a uniform maxmin value for player~\(i\).
Then, for every \(w\in \overline{W}\) and every \(\varepsilon>0\), there exist \(\delta>0\) and a strategy \(\sigma\in \Sigma^i\) such that \(\sigma\) is a robust \(\varepsilon\)-maxmin strategy for the set \(B(w,\delta)\cap W\).
\end{lem}

\textbf{Outline of the proof.} Let \(v^i\) be the uniform maxmin value of player~\(i\) in \(\Gamma\), and fix \(\varepsilon>0\).
Choose auxiliary positive parameters \(\varepsilon_0>0\), \(\varepsilon_1>0\), and \(\varepsilon_2>0\) (to be fixed later) with
\(\varepsilon_0+4\varepsilon_1+\varepsilon_2\le \varepsilon/2\).
By the definition of the uniform maxmin value, there exist an integer \(n_\varepsilon\) and a strategy \(\sigma^i_\varepsilon\in\Sigma^i\) such that, for every state \(z\in S\), the following two properties hold:
\[
\gamma^i_n(\sigma^i_\varepsilon,\tau)[z]\ge v^i[z]-\varepsilon_0
\quad\text{for every }\tau\in\Sigma^{-i} \quad\text{and } n\ge n_\varepsilon,
\]
and, for every \(\rho^i\in\Sigma^i\) and \(w\in W\) with \(w_1<1/n_\varepsilon\), there exists \(\tau\in\Sigma^{-i}\) such that
\[
\gamma^i_w(\rho^i,\tau)[z]\le v^i[z]+\varepsilon_0 .
\]
Thus, \(\sigma^i_\varepsilon\) is \(\varepsilon_0\)-maxmin in each \(n\)-stage game \(\Gamma_n\) for all \(n\ge n_\varepsilon\).

For any \(w\in\overline{W}\), fix a cutoff stage \(t_\varepsilon=t_\varepsilon(w)\) and define an auxiliary \(t_\varepsilon\)-stage game \(\Gamma_\ast=\Gamma_\ast(t_\varepsilon)\) with the same strategy spaces and transition kernel as \(\Gamma\), but with the payoff to player~\(i\) given by
\[
\sum_{t=1}^{t_\varepsilon-1} w_t\,g^i_t \;+\; \beta(w)\,v^i\bigl[z_{t_\varepsilon}\bigr],
\qquad
\text{where } \beta(w):=1-\sum_{t=1}^{\infty} w_t.
\]
(Here \(z_{t_\varepsilon}\) is the random state at stage \(t_\varepsilon\).)
Let \(\sigma_\ast\in\Sigma^i\) be an \(\varepsilon_2\)-maxmin strategy for player~\(i\) in \(\Gamma_\ast\).

We then construct a strategy for player~\(i\) that follows \(\sigma_\ast\) up to stage \(t_\varepsilon - 1\), and from stage \(t_\varepsilon\) onward switches to \(\sigma_\varepsilon\). The parameters \(\varepsilon_0\), \(t_\varepsilon\), and \(\varepsilon_2\) are chosen so that, for a sufficiently small \(\delta > 0\), the resulting strategy is a robust \(\varepsilon\)-maxmin strategy for the set of impatient valuations \(u\) with \(d(u, w) < \delta\), and for every strategy \(\rho\) of player~\(i\) and an impatient valuation \(u\) with \(d(u, w) < \delta\), there exists a strategy profile \(\tau\) of the other players so that \(\gamma^i_u(\rho,\tau)\le v^i(u)+\varepsilon\).

\begin{proof}
Suppose \(\Gamma\) has a uniform maxmin value \(v^i\) for player~\(i\). Fix \(w\in\overline{W}\) and \(\varepsilon>0\).
To prove the lemma, it suffices to construct \(\sigma\in\Sigma^i\) and \(\delta>0\) such that \(\sigma\) is a robust \(\varepsilon\)-maxmin strategy of player~\(i\) for the set
\[
B(w,\delta)\cap W \;=\; \{\,u\in W : d(u,w)<\delta\,\}.
\]

Choose auxiliary positive constants \( \varepsilon_0\), \(\varepsilon_1\), and \(\varepsilon_2 \) such that
\(\varepsilon_0+4\varepsilon_1+\varepsilon_2\le \varepsilon/2\).
One possible choice is \(\varepsilon_0=\varepsilon_1=\varepsilon_2=\varepsilon/12\), but we retain distinct symbols to reflect the structure of the error estimates below.

\smallskip

\noindent
\textbf{Definition of \( \beta \).}
Let \( \beta:=\beta(w)=1-\sum_{t=1}^\infty w_t \).
Since \(w\in\overline W\), we have \(0\le \sum_{t=1}^\infty w_t\le 1\); hence \(1\ge \beta\ge 0\).
We distinguish between two cases: \( \beta:=\beta(w)=0\) and
\( \beta:=\beta(w)>0\). The case \(\beta=0\), namely the case in which \(w\) has total mass one, is easy; nonetheless, for completeness, we present it.

\smallskip

\noindent
\textbf{Definition of the game \( \Gamma_\ast \).}
For an integer \(t_\varepsilon\) to be chosen below, let
\[
\Gamma_\ast=\Gamma_\ast(t_\varepsilon)
\]
be the strategic-form game with the same strategy spaces and transition kernel as
\(\Gamma\), and with player \(i\)'s payoff modified to
\[
\gamma^i_\ast(\rho)[z]
:=
E^z_\rho
\left[
\sum_{1\le t<t_\varepsilon} w_t g^i_t
+
\beta v^i[z_{t_\varepsilon}]
\right],
\]
where \(z_{t_\varepsilon}\) denotes the (random) state at stage \(t_\varepsilon\) and \(v^i\) is the uniform maxmin value of player~\(i\).
(We do not use the payoffs of the other players in \( \Gamma_\ast \); they may be defined arbitrarily as bounded functions.)

\smallskip
By standard arguments, \(\Gamma_\ast\) admits a maxmin value for player~\(i\): in the finite case, \(\Gamma_\ast\) is a finite normal-form game (von Neumann's minimax theorem); in the general case, \(\Gamma_\ast\) is a finite-horizon stochastic game, so backward induction yields a maxmin value. Hence there exists a value vector \(v_\ast^{\,i}\in\mathbb{R}^S\) and, for every \(\varepsilon_2>0\), an \(\varepsilon_2\)-maxmin strategy \(\sigma^i_\ast\in\Sigma^i\) in \(\Gamma_\ast\), i.e.
\[
\gamma_\ast^{\,i}(\sigma^i_\ast,\tau)[z]\ge v_\ast^{\,i}[z]-\varepsilon_2
\quad\text{for all }\tau\in\Sigma^{-i},\ z\in S.
\]
Moreover, for every strategy \(\hat{\sigma}^i_\ast\in\Sigma^i\), there exists
\(\tau\in\Sigma^{-i}\) such that
\[
\gamma_\ast^{\,i}(\hat{\sigma}^i_\ast,\tau)[z]\le v_\ast^{\,i}[z]+\varepsilon_2
\quad\text{for all }z\in S.
\]

\smallskip

\noindent
\textbf{Definition of \(v_\ast^{\,i}\) and \(\sigma_\ast\).}
Let \(\varepsilon_2>0\) be as fixed above. Whenever \(t_\varepsilon\) has been fixed, let \(v^i_\ast=v^i_\ast(t_\varepsilon)\) be the maxmin value vector of player~\(i\) in \(\Gamma_\ast(t_\varepsilon)\), and choose an \(\varepsilon_2\)-maxmin strategy \(\sigma_\ast=\sigma_\ast(t_\varepsilon)\in\Sigma^i\) for player~\(i\) in \(\Gamma_\ast(t_\varepsilon)\).\footnote{In the finite case, \(\sigma_\ast\) may be chosen maxmin optimal.}
Then
\begin{subequations}\label{eqn:maxmin_Gamma_star}
\begin{equation}\label{eqn:maxminsigmaitauminusi}
\gamma_\ast^{\,i}(\sigma_\ast,\tau)[z]\ge v_\ast^{\,i}[z]-\varepsilon_2
\quad\text{for all }\tau\in\Sigma^{-i},\ z\in S,
\end{equation}
and for every strategy \(\hat{\sigma}_\ast\in\Sigma^i\), there exists
\(\hat{\tau}_\ast\in\Sigma^{-i}\) such that
\begin{equation}\label{eqn:countermaxminsigmaitauminusi}
\gamma_\ast^{\,i}(\hat{\sigma}_\ast,\hat{\tau}_\ast)[z]\le v_\ast^{\,i}[z]+\varepsilon_2
\quad\text{for all }z\in S.
\end{equation}
\end{subequations}

If \(\beta=0\), \(t_\varepsilon\) is large enough so that \(\|g^i\|\sum_{t=t_\varepsilon}^\infty w_t<\varepsilon_1\), and \(\delta=\delta(w,t_\varepsilon)>0\) is small enough so that
\[
\|g^i\|
\left(
\sum_{t=t_\varepsilon}^{\infty}w_t+2t_\varepsilon\delta
\right)
<\varepsilon_1,
\]
and \(\sigma\) plays according to \(\sigma_\ast\) through stage \(t_\varepsilon-1\), then, for every \(u\in B(w,\delta)\cap W\),
\[
\sum_{1\leq t<t_\varepsilon} |g^i_t(u_t-w_t)|+\sum_{t=t_\varepsilon}^\infty |g^i_t u_t|<\varepsilon_1,
\]
\begin{equation}\label{eqn:maxminsigmaitauminusibetaequalzero}
\gamma_u^{\,i}(\sigma,\tau)[z]\ge v_\ast^{\,i}[z]-\varepsilon_2-\varepsilon_1\ge v_\ast^{\,i}[z]-\varepsilon/2
\quad\text{for all }\tau\in\Sigma^{-i},\ z\in S,
\end{equation}
and for every strategy \(\rho\) of player~\(i\) there is a strategy profile \(\tau\) of the other players such that

\begin{equation}\label{eqn:countertaumaxmin_betaequalzero}
\begin{aligned}
\gamma^i_u(\rho,\tau)[z]
&\le
E^z_{\rho,\tau}
\left[
\sum_{1\le t<t_\varepsilon} w_t g^i_t
\right]
+
\sum_{1\le t<t_\varepsilon}|g^i_t(u_t-w_t)|
+
\sum_{t=t_\varepsilon}^{\infty}|g^i_tu_t|  \\
&\le
v^i_\ast[z]+\varepsilon_2+\varepsilon_1
\le
v^i_\ast[z]+\varepsilon/2
\qquad\text{for all } z\in S .
\end{aligned}
\end{equation}
Thus, for every \(u\in B(w,\delta)\cap W\) and every \(z\in S\),
\[
v^i_\ast[z]-\varepsilon/2
\le
v^i(u)[z]
\le
v^i_\ast[z]+\varepsilon/2.
\]
Together with the lower bound in \eqref{eqn:maxminsigmaitauminusibetaequalzero}, this implies that \(\sigma\) is a robust \(\varepsilon\)-maxmin strategy for \(B(w,\delta)\cap W\).

We assume, hereunder, that \(\beta>0\). In this case, the cutoff \(t_\varepsilon\) appearing in the auxiliary game \(\Gamma_\ast=\Gamma_\ast(t_\varepsilon)\) will be chosen
after \(n_\varepsilon\) has been fixed, and we will then impose the additional constraints on \(t_\varepsilon\) and \(\delta\) that are needed below.

\smallskip

\noindent
\textbf{Definition of \( \sigma_\varepsilon \) and \( n_\varepsilon \).}
By the assumption that player~\(i\) has a uniform maxmin value \(v^i\), there exist a strategy \( \sigma_\varepsilon \in \Sigma^i \) and an integer \( n_\varepsilon \) such that
\begin{equation}\label{eqn:uniformmaxminvalue}
\gamma^i_n(\sigma_\varepsilon,\tau)[z] \ge v^i[z]-\varepsilon_0
\quad \text{for all } n \ge n_\varepsilon,\ \tau \in \Sigma^{-i},\ z \in S.
\end{equation}

We now choose \(t_\varepsilon\) sufficiently large, and then \(\delta>0\) sufficiently small.

\smallskip

\noindent
\textbf{Definition of \( t_\varepsilon \) and \(\delta\).}
Note that since \(w\in\overline W\),
\[
(n_\varepsilon+s)\,w_s+\sum_{t=s}^\infty w_t \longrightarrow 0 \quad \text{as } s\to\infty,
\]
because \(w_t\) is nonincreasing and \(\sum_{t=1}^\infty w_t<\infty\).

Let \(\theta>0\) be the threshold supplied by Lemma~\ref{lem:cesaro_patient}
with \(M=\|g^i\|\), \(N=n_\varepsilon+1\), \(\eta=\varepsilon_0\), and tolerance \(\varepsilon_0+2\varepsilon_1\).
We choose \(t_\varepsilon\) large enough so that
\[
\|g^i\|\sum_{t=t_\varepsilon}^{\infty}w_t
+
\|g^i\|(n_\varepsilon+t_\varepsilon)w_{t_\varepsilon}
<\varepsilon_1,\quad n_\varepsilon w_{t_\varepsilon}<\beta,
\quad \text{and }
w_{t_\varepsilon}<\theta.
\]

Then choose \(\delta>0\) small enough so that
\begin{equation}\label{eqn:condition_on_delta}
\|g^i\|\sum_{t=t_\varepsilon}^{\infty}w_t
+
\|g^i\|(n_\varepsilon+t_\varepsilon)(w_{t_\varepsilon}+\delta)
<\varepsilon_1,\quad n_\varepsilon(w_{t_\varepsilon}+\delta)<\beta-t_\varepsilon\delta,
\quad \text{and } w_{t_\varepsilon}+\delta<\theta.
\end{equation}

It follows that for \(u\in B(w,\delta)\cap W\)
\[
u_{t_\varepsilon}\le w_{t_\varepsilon}+\delta<\theta \quad \text{and } \sum_{t=t_\varepsilon}^{\infty}u_t
=
1-\sum_{t<t_\varepsilon}u_t
\ge
1-\sum_{t<t_\varepsilon}w_t-t_\varepsilon\delta
\ge
\beta-t_\varepsilon\delta .
\]
Let \(u\in B(w,\delta)\cap W\), and set
\[
s_u:=\sum_{t=t_\varepsilon}^{\infty}u_t .
\]
By the choice of \(t_\varepsilon\) and \(\delta\), \(s_u>0\). Define the normalized tail valuation
\[
\widehat u_r:=\frac{u_{t_\varepsilon+r-1}}{s_u},
\qquad r\ge1.
\]
Then \(\widehat u\in W\), and
\begin{equation}\label{eqn:normalized_tail_patient}%
\widehat u_1
=
\frac{u_{t_\varepsilon}}{\sum_{t=t_\varepsilon}^{\infty}u_t}
\le
\frac{w_{t_\varepsilon}+\delta}{\beta-t_\varepsilon\delta}
<
\frac1{n_\varepsilon}.
\end{equation}

\smallskip

\noindent
\textbf{Definition of \( \sigma \).}
Define \(\sigma\in\Sigma^i\) to play according to \(\sigma_\ast\) through stage \(t_\varepsilon-1\), and from stage \(t_\varepsilon\) onward to follow \(\sigma_\varepsilon\) in the continuation game whose initial state is \(z_{t_\varepsilon}\).
Formally,
\begin{itemize}
  \item for \(1\le t<t_\varepsilon\), \(\ \sigma(z_1,\ldots,z_t)=\sigma_\ast(z_1,\ldots,z_t)\);
\item for \(t\ge t_\varepsilon\), \(\ \sigma(z_1,\ldots,z_t)=\sigma_\varepsilon(z_{t_\varepsilon},\ldots,z_t)\) (i.e., it depends only on the suffix starting at \(z_{t_\varepsilon}\)).
\end{itemize}

\smallskip
\smallskip

\noindent
\textbf{Proof that \( \sigma \) is a robust \( \varepsilon \)-maxmin strategy of player~\( i \) for the set of impatient valuations in \( B(w, \delta) \).}
Let \( u \in B(w, \delta) \) be an impatient valuation, and let \( \tau \in \Sigma^{-i} \).
By the definition of \(\sigma\), for every \(\rho\in\Sigma^{-i}\), the distributions
\(P^z_{\sigma_\ast,\rho}\) and \(P^z_{\sigma,\rho}\)
on \((z_1,a_1,\ldots,z_{t_\varepsilon})\) coincide:
\begin{equation}\label{eqn:prefix_distribution_coincidence}
P^z_{\sigma_\ast,\rho}=P^z_{\sigma,\rho}
\quad\text{on }(z_1,a_1,\ldots,z_{t_\varepsilon}).
\end{equation}

From \eqref{eqn:prefix_distribution_coincidence} and \eqref{eqn:maxminsigmaitauminusi}, together with the pointwise bounds
\[
w_t g^i_t + \|g^i\|\,|u_t-w_t|\ge u_t g^i_t \ge w_t g^i_t - \|g^i\|\,|u_t-w_t|,
\]
the inequality  \( |u_t-w_t|<\delta\), and using \eqref{eqn:condition_on_delta} (which implies \( \|g^i\|\,t_\varepsilon \delta < \varepsilon_1 \)), we obtain
\begin{subequations}\label{eqn:maxmin_tau}
\begin{align}
E^z_{\sigma, \tau} \!\left[ \sum_{1 \le t < t_\varepsilon} u_t g^i_t \right]
&\ge E^z_{\sigma_\ast, \tau} \!\left[ \sum_{1 \le t < t_\varepsilon} w_t g^i_t \right]
    - \|g^i\| \sum_{1 \le t < t_\varepsilon} |w_t - u_t| \nonumber \\
&\ge v_\ast^i[z] - \varepsilon_2 - \beta\, E^z_{\sigma, \tau} v^i[z_{t_\varepsilon}] - \varepsilon_1,
\label{eqn:maxmintau1new}
\intertext{and for every \(\hat{\sigma}\in \Sigma^i\) there exists \(\hat{\tau}_\ast\in\Sigma^{-i}\) such that}
E^z_{\hat{\sigma}, \hat{\tau}_\ast} \!\left[ \sum_{1 \le t < t_\varepsilon} u_t g^i_t \right]
&\le v_\ast^i[z] + \varepsilon_2 - \beta\, E^z_{\hat{\sigma}, \hat{\tau}_\ast} v^i[z_{t_\varepsilon}] + \varepsilon_1.
\label{eqn:countermaxmintau1new}
\end{align}
\end{subequations}

By \eqref{eqn:uniformmaxminvalue} and the definition of \(\sigma\), for every
\(t\ge t_\varepsilon+n_\varepsilon\),
\begin{equation}\label{eqn:maxminsigmaFt4}
E_{\sigma,\tau}
\left(
g^i_{t_\varepsilon}+\cdots+g^i_t
\,\middle|\,
\mathcal H_{t_\varepsilon}
\right)
\ge
(t-t_\varepsilon+1)\bigl(v^i[z_{t_\varepsilon}]-\varepsilon_0\bigr).
\end{equation}

For every \(\hat\sigma\in\Sigma^i\), apply the second condition of
Definition~\ref{def:uniform_maxmin} to the normalized tail valuation
\(\widehat u\). Since \(\widehat u_1<1/n_\varepsilon\), there exists a continuation
counter-strategy \(\hat\tau^{\,\mathrm{tail}}\in\Sigma^{-i}\) such that
\[
E_{\hat\sigma,\hat\tau^{\,\mathrm{tail}}}
\left[
\sum_{t=t_\varepsilon}^{\infty}u_t g^i_t
\,\middle|\,
\mathcal H_{t_\varepsilon}
\right]
\le
\sum_{t=t_\varepsilon}^{\infty}u_t
\bigl(v^i[z_{t_\varepsilon}]+\varepsilon_0\bigr).
\]
Since \(\sum_{t=t_\varepsilon}^{\infty}u_t\le 1\), it follows that
\begin{equation}\label{eqn:countermaxminsigmaFt4}
E_{\hat\sigma,\hat\tau^{\,\mathrm{tail}}}
\left[
\sum_{t=t_\varepsilon}^{\infty}u_t g^i_t
\,\middle|\,
\mathcal H_{t_\varepsilon}
\right]
\le
\sum_{t=t_\varepsilon}^{\infty}u_t v^i[z_{t_\varepsilon}]
+\varepsilon_0 .
\end{equation}

As \( |v^i[z_{t_\varepsilon}]| \le \|g^i\| \) and \( \delta>0 \) was chosen so that \(\|g^i\|\sum_{t=t_\varepsilon}^\infty w_t+ \|g^i\|\,t_\varepsilon \delta < \varepsilon_1 \), we have
\begin{subequations}
\begin{align}
\sum_{t = t_\varepsilon}^\infty u_t\, v^i[z_{t_\varepsilon}]
&= \Bigl(1 - \sum_{t = 1}^{t_\varepsilon - 1} u_t \Bigr) v^i[z_{t_\varepsilon}]
\ge \Bigl(1 - \sum_{t = 1}^{t_\varepsilon - 1} w_t \Bigr) v^i[z_{t_\varepsilon}] - \|g^i\|t_\varepsilon\,\delta \nonumber\\
&=\beta v^i[z_{t_\varepsilon}]+\sum_{t=t_\varepsilon}^\infty w_tv^i[z_{t_\varepsilon}]- \|g^i\|t_\varepsilon\,\delta
\ge \beta\, v^i[z_{t_\varepsilon}] - \varepsilon_1.
\label{eqn:lower_u_t_v_i}
\intertext{and}
\sum_{t = t_\varepsilon}^\infty u_t\, v^i[z_{t_\varepsilon}]
&= \Bigl(1 - \sum_{t = 1}^{t_\varepsilon - 1} u_t \Bigr) v^i[z_{t_\varepsilon}]
\le \Bigl(1 - \sum_{t = 1}^{t_\varepsilon - 1} w_t \Bigr) v^i[z_{t_\varepsilon}] + \|g^i\|t_\varepsilon\,\delta \nonumber\\
&=\beta v^i[z_{t_\varepsilon}]+\sum_{t=t_\varepsilon}^\infty w_tv^i[z_{t_\varepsilon}]+ \|g^i\|t_\varepsilon\,\delta
\le \beta\, v^i[z_{t_\varepsilon}] + \varepsilon_1.
\label{eqn:upper_u_t_v_i}
\end{align}
\end{subequations}

\smallskip

Since \(u_{t_\varepsilon}<\theta\), Lemma~\ref{lem:cesaro_patient} applies to the nonincreasing sequence \(\lambda_r:=u_{t_\varepsilon+r-1}\).

Applying Lemma~\ref{lem:cesaro_patient} pathwise, after conditioning on \(\mathcal{H}_{t_\varepsilon}\), to the sequence
\[ x_r := E_{\sigma,\tau}\!\left( g^i_{t_\varepsilon+r-1} \,\middle|\, \mathcal{H}_{t_\varepsilon} \right), \qquad r\ge 1, \]
with \[ a:=v^i[z_{t_\varepsilon}],\qquad \eta:=\varepsilon_0,\qquad N:=n_\varepsilon+1,\qquad \lambda_r:=u_{t_\varepsilon+r-1}, \] and with tolerance \(\varepsilon_0+2\varepsilon_1\), the bound \eqref{eqn:maxminsigmaFt4} implies \begin{equation}\label{eqn:maxminsigmaFt5} E_{\sigma,\tau}\!\left( \sum_{t=t_\varepsilon}^\infty u_t g^i_t \,\middle|\, \mathcal{H}_{t_\varepsilon} \right) \ge \sum_{t=t_\varepsilon}^\infty u_t\,v^i[z_{t_\varepsilon}] -\varepsilon_0-2\varepsilon_1 . \end{equation}

Therefore, by \eqref{eqn:maxminsigmaFt5}, \eqref{eqn:lower_u_t_v_i}, and \eqref{eqn:prefix_distribution_coincidence}
\begin{equation}\label{eqn:maxmintau2new}
E_{\sigma, \tau} \sum_{t = t_\varepsilon}^\infty u_t g^i_t
\;\ge\; \beta\, E_{\sigma_\ast, \tau} v^i[z_{t_\varepsilon}] - \varepsilon_0 - 3\varepsilon_1
\;=\; \beta\, E_{\sigma, \tau} v^i[z_{t_\varepsilon}] - \varepsilon_0 - 3\varepsilon_1 .
\end{equation}

Adding \eqref{eqn:maxmintau1new} and \eqref{eqn:maxmintau2new} yields
\[
\gamma^i_u(\sigma, \tau)[z]
= E^z_{\sigma, \tau} \sum_{t = 1}^\infty u_t g^i_t
\ge v_\ast^i[z] - \varepsilon_0 - \varepsilon_2 - 4\varepsilon_1 .
\]
Since \( \varepsilon_0 + 4\varepsilon_1 + \varepsilon_2 \le \varepsilon/2 \), it follows that
\[
\gamma^i_u(\sigma, \tau)[z] \ge v_\ast^i[z] - \varepsilon/2,
\]

We now prove the corresponding upper bound. Fix \(\hat\sigma\in\Sigma^i\). Let \(\hat\tau_\ast\) be the counter-strategy supplied by \eqref{eqn:countermaxmintau1new} in the auxiliary game
\(\Gamma_\ast\). For each history at stage \(t_\varepsilon\), let
\(\hat\tau^{\,\mathrm{tail}}\) be the continuation counter-strategy supplied by the second condition of Definition~\ref{def:uniform_maxmin} for the normalized tail valuation
\(\widehat u\). Define \(\hat\tau\in\Sigma^{-i}\) to follow \(\hat\tau_\ast\) in stages
\(1\le t<t_\varepsilon\), and from stage \(t_\varepsilon\) onward to use
\(\hat\tau^{\,\mathrm{tail}}\).

Then \eqref{eqn:countermaxmintau1new} and \eqref{eqn:countermaxminsigmaFt4} hold with the spliced strategy \(\hat\tau\).
Using these inequalities together with \eqref{eqn:upper_u_t_v_i}, we obtain
\[
\begin{aligned}
\gamma^i_u(\hat\sigma,\hat\tau)[z]
&\le
v^i_\ast[z]+\varepsilon_2
-\beta E^z_{\hat\sigma,\hat\tau}v^i[z_{t_\varepsilon}]
+\varepsilon_1  \\
&\quad+
\beta E^z_{\hat\sigma,\hat\tau}v^i[z_{t_\varepsilon}]
+\varepsilon_1+\varepsilon_0  \\
&=
v^i_\ast[z]+\varepsilon_0+\varepsilon_2+2\varepsilon_1
\le
v^i_\ast[z]+\varepsilon/2 .
\end{aligned}
\]
Thus, for every \(u\in B(w,\delta)\cap W\) and every \(z\in S\),
\[
v_\ast^i[z]-\varepsilon/2 \le v^i(u)[z]\le v_\ast^i[z]+\varepsilon/2 .
\]
Combining this with the lower bound
\[
\gamma_u^i(\sigma,\tau)[z]\ge v_\ast^i[z]-\varepsilon/2
\]
shows that, for all \(u,u'\in B(w,\delta)\cap W\), all \(\tau\in\Sigma^{-i}\), and all \(z\in S\),
\[
\gamma_u^i(\sigma,\tau)[z]\ge v^i(u')[z]-\varepsilon .
\]
Hence \(\sigma\) is a robust \(\varepsilon\)-maxmin strategy of player~\(i\) for \(B(w,\delta)\cap W\).
\end{proof}

\begin{lem}\label{lem:accumulationofUn}
Let \( U_n \subset W \) be a sequence of nonempty sets whose diameters converge to zero. Then there exists \( w \in \overline{W} \) such that every neighborhood of \( w \) contains some \( U_n \).
\end{lem}

\begin{proof}
We first show that \( \overline{W} \) is sequentially compact. Every sequence in \( \overline{W} \) is bounded, and since each element \( w = (w_t)_{t \geq 1} \in \overline{W} \) satisfies \( 0 \leq w_t \leq 1/t \), the diagonal argument yields a subsequence that converges pointwise. Moreover, since pointwise limits of nonincreasing sequences with \( \sum_{t=1}^\infty w_t \leq 1 \) also lie in \( \overline{W} \), the limit is valid. Finally, because \( w_t \leq 1/t \) implies \( w_t \to 0 \), the convergence is uniform, and hence pointwise convergence implies convergence in the supremum metric \( d \). Therefore, \( \overline{W} \) is sequentially compact.

Now let \( w(n) \in U_n \). By sequential compactness, there exists a subsequence \( w(n_k) \to w \in \overline{W} \).
Since \( \operatorname{diam}(U_{n_k}) \to 0 \) as \( k \to \infty \), it follows that for every neighborhood \( B(w, \delta) \), there exists \( N \) such that \( U_{n_k} \subset B(w, \delta) \) for all \( k \geq N \). Thus, every neighborhood of \( w \) contains some \( U_n \).
\end{proof}

Assume that condition (1.a) holds; that is, \(\Gamma\) admits a uniform maxmin value for player \(i\).
By Lemma~\ref{lem:uniformmaxminatwinoverlineW}, \(\Gamma\) has robust \(\varepsilon\)-maxmin strategies for player \(i\) at every
\(w\in\overline{W}\). Suppose condition (1.b) does not hold.
Then there exist \(\varepsilon>0\) and a sequence \((U_n)\subset W\) of nonempty subsets whose diameters converge to zero such that, for every \(n\), no robust \(\varepsilon\)-maxmin strategy exists for
\(U_n\).
By Lemma~\ref{lem:accumulationofUn},
there exists \(w\in \overline{W}\)
such that every neighborhood of \(w\) contains some \(U_n\).
By local robustness at \(w\)
 (Lemma~\ref{lem:uniformmaxminatwinoverlineW})
there exist \(\delta_w>0\) and a strategy that is robust \(\varepsilon\)-maxmin on \(B(w,\delta_w)\cap W\); since every
neighborhood of \(w\) contains some \(U_n\),
some \(U_n\subseteq B(w,\delta_w)\).
That \(U_n\) therefore
admits a robust \(\varepsilon\)-maxmin strategy, contradicting the choice of \((U_n)\). Hence condition (1.b) holds.
\qed

\smallskip
Assume that condition (1.b) holds, and fix \( \varepsilon > 0 \). Then there exists \( \delta > 0 \) such that, for every \( w \in \overline{W} \), player~\( i \) has a robust \( \varepsilon \)-maxmin strategy for the set \( B(w, 2\delta) \cap W \).

Since the metric space \( \overline{W} \) is compact, and since \( \bigcup_{w \in \overline{W}} B(w, \delta) = \overline{W} \), there exists a finite list \( w(1), \ldots, w(k) \in \overline{W} \) such that the corresponding balls \( B(w(j), \delta) \) cover \( \overline{W} \).

For each \( j \), let \( \sigma_j \) be a robust \( \varepsilon \)-maxmin strategy of player~\( i \) for \( B(w(j), 2\delta) \cap W \), and define \( U_j := B(w(j), \delta) \cap W \). Then \( \sigma_j \) is a robust \( \varepsilon \)-maxmin strategy for \( U_j^\delta \), and \( \bigcup_j U_j = W \).

This completes the proof that condition (1.b) implies condition (1.c). \qed

\smallskip
Assume that condition (1.c) holds, and fix \( \varepsilon > 0 \). Then there exists \( \delta > 0 \) and a finite list \( (\sigma_j, U_j) \), \( 1 \leq j \leq k \), where each \( \sigma_j \) is a strategy of player~\( i \), and each \( U_j \subset W \), such that \( \bigcup_j U_j = W \), and \( \sigma_j \) is a robust \( \varepsilon \)-maxmin strategy for player~\( i \) on \( U_j^{2\delta} \).

For each \( j \), let  \( f_j : \overline{W} \to \mathbb{R}_+ \) be a continuous function such that \( f_j(w) > 0 \)  if and only if \( w \in B(U_j, \delta) \), where
\[
B(U_j, \delta) := \{ w \in \overline{W} : d(w, U_j) < \delta \}.
\]
For example, one may define \( f_j(w) := \mathrm{dist}(w, \overline{W} \setminus B(U_j, \delta)) \), which is continuous, strictly positive on \( B(U_j, \delta) \), and vanishes outside \( B(U_j, \delta) \). Since the sets \( B(U_j, \delta) \) cover \( \overline{W} \), the sum \( \sum_j f_j(w)  \) is strictly positive for all \( w \in \overline{W} \).

Define weights
\[
\alpha_j(w) := \frac{f_j(w)}{\sum_{l=1}^k f_l(w)},
\]
which form a continuous partition of unity. Then, the mapping
\[
\overline{W} \ni w \mapsto \sigma_w := \sum_{j=1}^k \alpha_j(w) \sigma_j
\]
is uniformly continuous, as it is at each point \(w\in \overline{W}\) a convex combination of a finite collection of fixed strategies, with weights that vary uniformly continuously in \( w \) over the compact domain \( \overline{W} \).

Now fix \( w \in W \). For each \( j \) with \( f_j(w) > 0 \), we have \( d(w, U_j) < \delta \), and hence \( B(w, \delta) \subseteq B(U_j, 2\delta) \). Therefore, each \( \sigma_j \) with positive weight (in the mixture) is a robust \( \varepsilon \)-maxmin strategy for \( B(w, \delta) \cap W \). Hence, the convex combination \( \sigma_w \) is also a robust \( \varepsilon \)-maxmin strategy for \( B(w, \delta) \cap W \).

It follows that \( \sigma_w \) is a robust \( \varepsilon \)-maxmin strategy of player~\( i \) for \( \{w\}^\delta = B(w, \delta) \cap W \), and the mapping \( w \mapsto \sigma_w \) is uniformly continuous.

This completes the proof that condition~(1.c) implies condition~(1.d). \qed

\smallskip
Assume that condition \((1.d)\) holds. We first show that \(v^i_n\) converges
in the sup norm over \(S\). Let \(\eta>0\). By condition \((1.d)\), there exist
\(\delta>0\) and, for every \(w\in W\), a strategy \(\sigma_w\) that is a robust
\(\eta\)-maxmin strategy for \(B(w,\delta)\cap W\). Choose \(N\) so large that
\(1/N<\delta\). Then, for all \(m,n\ge N\), the Cesaro valuations \(u_m\) and
\(u_n\) both belong to \(B(u_m,\delta)\cap W\). Therefore, the value variation
under robust maxmin strategies gives
\[
\|v^i_m-v^i_n\|_\infty\le \eta .
\]
Thus \((v^i_n)_n\) is a Cauchy sequence in the Banach space \(\ell^\infty(S)\),
and hence converges in sup norm to some bounded function \(v^i:S\to\mathbb R\).

We now verify that \(v^i\) satisfies the two conditions of Definition~\ref{def:uniform_maxmin}.
Fix \(\varepsilon>0\), and set \(\eta:=\varepsilon/3\). Apply condition \((1.d)\)
with tolerance \(\eta\), and let \(\delta>0\) and \(w\mapsto\sigma_w\) be as above.
Choose \(k\) large enough so that \(1/k<\delta\) and
\[
\|v^i_n-v^i\|_\infty<\eta
\quad\text{for all } n\ge k.
\]
Set \(n_\varepsilon:=k\) and \(\sigma^i_\varepsilon:=\sigma_{u_k}\).

For every \(n\ge n_\varepsilon\), we have \(u_n\in B(u_k,\delta)\cap W\). Hence,
for every \(\tau\in\Sigma^{-i}\) and every \(z\in S\),
\[
\gamma^i_n(\sigma^i_\varepsilon,\tau)[z]
\ge v^i_n[z]-\eta
\ge v^i[z]-2\eta
\ge v^i[z]-\varepsilon .
\]
This establishes the first condition of Definition~\ref{def:uniform_maxmin}.

We next prove the second condition. Let \(w\in W\) satisfy \(w_1<1/n_\varepsilon\).
Since \(1/n_\varepsilon<\delta\), we have
\[
\|w-u_{n_\varepsilon}\|_\infty<1/n_\varepsilon<\delta,
\]
and therefore \(w\in B(u_{n_\varepsilon},\delta)\cap W\). The value variation
under robust maxmin strategies gives
\[
\|v^i(w)-v^i_{n_\varepsilon}\|_\infty\le \eta .
\]
Therefore, for every strategy \(\rho^i\in\Sigma^i\), there exists
\(\tau\in\Sigma^{-i}\) such that
\[
\gamma^i_w(\rho^i,\tau)[z]
\le v^i(w)[z]+\eta
\le v^i_{n_\varepsilon}[z]+2\eta
\le v^i[z]+3\eta
= v^i[z]+\varepsilon
\]
for every \(z\in S\). This establishes the second condition of
Definition~\ref{def:uniform_maxmin}. Hence \(v^i\) is the uniform maxmin value of player~\(i\).

This proves that condition (1.d) implies condition (1.a), and completes the proof of the equivalence of the stated conditions in Theorem~\ref{thm:robustoptimalstrategy}.
\qed

\subsection{Proof of Theorem \ref{thm:uniformcontinuityofv}}
Every finite stochastic game \(\Gamma\) has a uniform maxmin value for each player~\(i\); see \cite{ney03minmax}, building on \cite{merney81}. Hence, condition~(1.a) of Theorem~\ref{thm:robustoptimalstrategy} holds. By the same theorem,  condition~(1.b) also holds. That is, for every \( \varepsilon > 0 \), there exists \( \delta > 0 \) such that for any subset \( U \subset W \) with \( \operatorname{diam}(U) < \delta \), player~\( i \) has a robust \( \varepsilon \)-maxmin strategy \( \sigma^i \) for \( U \).
This implies that
\[
v^i(u)[z] \geq \inf_{\sigma^{-i} \in \Sigma^{-i}} \gamma^i_u(\sigma^i, \sigma^{-i}) \geq v^i(w)[z] - \varepsilon
\quad \text{for all } u, w \in U, \; z \in S,
\]
which proves that the map \( W^I \ni w \mapsto v^i(w) \) is uniformly continuous.
\qed

\begin{remark}
In a finite stochastic game, each player~\( i \) has a Markov maxmin strategy in \( \Gamma_u \). Thus, the maxmin value \( v^i(u) \) is unaffected by whether or not player~\( i \) observes the other players' past actions. By contrast, the existence of a uniform \( \varepsilon \)-maxmin strategy generally requires that player~\( i \) observe the opponents’ past actions. Our proof of Theorem~\ref{thm:uniformcontinuityofv} relies on the existence of such uniform \( \varepsilon \)-maxmin strategies.

An alternative proof of the uniform continuity of the mapping \( u \mapsto v^i(u) \) avoids this reliance. It combines the compactness of \( \overline{W} \) with the existence of the limit of \( v^i(u) \) as \( u_1 \to 0 \), established in~\cite{zil2016,zil2018} without assuming the existence of a uniform maxmin. Specifically:
\begin{itemize}
    \item For every \( w \in \overline{W} \) and \( \varepsilon > 0 \), one shows that there exists \( \delta > 0 \) such that for all \( u, u' \in B(w, 2\delta) \cap W \), we have \( |v^i(u) - v^i(u')| < \varepsilon \).
    \item Then, using the compactness of \( \overline{W} \), one deduces that the variation of \( v^i(u) \) over any set of impatient valuations of diameter less than \( \delta \) is at most \( \varepsilon \), yielding uniform continuity.
\end{itemize}
\end{remark}

\subsection{Proof of Theorem \ref{thm:uniformcontinuousoptimalstrategies}}
A finite stochastic game \( \Gamma \) admits a uniform maxmin value for every player \( i \); see \cite{ney03minmax}, building on \cite{merney81}. Hence, condition~(1.a) of Theorem~\ref{thm:robustoptimalstrategy} holds. By the same theorem, this implies condition~(1.d), which asserts the existence of a uniformly continuous selection of \( \varepsilon \)-maxmin strategies.
\qed

The second part of this section focuses on the non-zero-sum case. The structure of the proof of Theorem~\ref{thm:robustequilibrium} parallels that of the adversary case, but crucial differences arise due to the nonconvexity of equilibrium sets.
The central idea is to show that the existence of uniform approximate equilibrium payoffs implies, for every \(\varepsilon>0\), the existence of a finite open cover \((U_j)_{j\in J}\) of the space of impatient valuation profiles \(W^I\), along with a corresponding list of strategy profiles \((\sigma_j)_{j\in J}\), such that each \(\sigma_j\) is a robust \(\varepsilon\)-equilibrium for all valuations in a neighborhood of its associated set \(U_j\).

The construction again relies on the compactness of the completion
\(\overline{W}^I\), which ensures the existence of a finite open cover and corresponding strategy profiles. However, unlike in the adversary case, we cannot exploit convexity to construct a globally continuous selection. Instead, we obtain approximate patchwise continuity: each strategy profile \(\sigma_j\) ensures robust equilibrium behavior over a local neighborhood, yielding a finite piecewise structure sufficient to characterize strategic robustness under preference uncertainty.

\subsection{Proof of Theorem  \ref{thm:robustequilibrium}}
As in the proof of Theorem~\ref{thm:robustoptimalstrategy}, we begin with the following intermediate result, which is of independent interest and establishes local robustness from the existence of uniform approximate equilibrium payoffs.

\begin{lem}\label{lem:uniformequilibriumatwinoverlineW}
Let \(\Gamma\) be a stochastic game that has uniform approximate equilibrium payoffs. Then, for every \(w \in \overline{W}^I\) and every \(\varepsilon > 0\), there exists \(\delta > 0\) and a strategy profile \(\sigma\) such that \(\sigma\) is a robust \(\varepsilon\)-equilibrium  for the set \(B(w, \delta) \cap W^I\).
\end{lem}

The lemma is the equilibrium counterpart of the maxmin Lemma~\ref{lem:uniformmaxminatwinoverlineW}. Its proof is analogous to the proof of Lemma~\ref{lem:uniformmaxminatwinoverlineW}.

\smallskip
\textbf{Outline of the proof.}
Assume that the stochastic game \(\Gamma\) has  uniform approximate equilibrium payoffs. Fix \(\varepsilon > 0\), define sufficiently small strictly positive constants \( \varepsilon_0\), \(\varepsilon_1\), and \( \varepsilon_2 \) such that \( \varepsilon_0 + 4\varepsilon_1 + \varepsilon_2 \leq \varepsilon/2 \).
By condition \((2.a^\ast)\), choose a bounded measurable payoff vector
\(v=(v^i[z])_{i\in I,z\in S}\), an integer \(n_\varepsilon\), and a strategy profile \(\sigma_\varepsilon\) such that, for every \(n\ge n_\varepsilon\), every \(i\in I\), every \(z\in S\), and every strategy \(\tau^i\) of player \(i\),
\[
\gamma^i_n(\sigma_\varepsilon)[z]\ge v^i[z]-\varepsilon_0
\quad\text{and}\quad
\gamma^i_n(\sigma_\varepsilon^{-i},\tau^i)[z]\le v^i[z]+\varepsilon_0 .
\]

By truncating \(v\) componentwise if necessary, we may assume without loss of generality that
\[
\|v^i\|_\infty\le \|g^i\|\qquad\text{for every }i\in I .
\]
The truncation preserves the two inequalities above and, in the infinite-state case, preserves measurability.

For each player \(i\), let \(\theta_i>0\) be the constant supplied by
Lemma~\ref{lem:cesaro_patient} with
\[
M=\|g^i\|,\qquad N=n_\varepsilon+1,\qquad
\eta=\varepsilon_0,
\]
and tolerance \(\varepsilon_0+2\varepsilon_1\).

Fix \(w \in \overline{W}^I\). Choose \(t_\varepsilon\) sufficiently large so that, for every \(i\in I\),
\begin{equation}\label{choosing-te}
\|g^i\|\sum_{t=t_\varepsilon}^{\infty} w^i_t
+
\|g^i\|(n_\varepsilon+t_\varepsilon)w^i_{t_\varepsilon}
<\varepsilon_1,
\qquad
w^i_{t_\varepsilon}<\theta_i/2 .
\end{equation}

Define the strategic-form game  \( \Gamma_\ast := (I, \Sigma, \gamma_\ast) \) in which each player's strategy set coincides with that of the original stochastic game, and player~\( i \)'s payoff is given by
\begin{equation}\label{eqn:defgammaistar}
\gamma^i_\ast(\sigma)[z] := E^z_\sigma \left( \sum_{1 \leq t < t_\varepsilon} w^i_t g^i_t + \beta_i v^i[z_{t_\varepsilon}] \right)
\quad \text{for all } \sigma \in \Sigma, \; z \in S,
\end{equation}
where \( \beta_i :=  1 - \sum_{t=1}^\infty w^i_t \).

Standard equilibrium arguments apply to the game \( \Gamma_\ast \), which is either, in the finite case, a game with finitely many pure strategies, or, in the general stochastic game case, a finite-horizon game in which, for every \( \varepsilon > 0 \), an \( \varepsilon \)-equilibrium exists. Therefore, the game \( \Gamma_\ast \) admits an \( \varepsilon_2\)-equilibrium. Let \( \sigma_\ast \) be an \( \varepsilon_2 \)-equilibrium strategy profile of \( \Gamma_\ast \), and let \( v_\ast = (v^i_\ast)_{i \in I} \) denote the corresponding payoff vector. That is,
\begin{equation}\label{eqn:sigmaFt3new}
\gamma^i_\ast(\sigma_\ast)[z] = v^i_\ast[z] \quad \text{for all } i \in I, \; z \in S,
\end{equation}
and for every player \( i \in I \) and every strategy \( \tau^i \) of player \( i \),
\begin{equation}\label{eqn:sigmaFt3bnew}
\gamma^i_\ast(\sigma_\ast^{-i}, \tau^i)[z] \leq v^i_\ast[z] + \varepsilon_2 \quad \text{for all }  z \in S.
\end{equation}

The strategy profile \( \sigma \) plays according to \( \sigma_\ast \) up to stage \( t_\varepsilon \), and then switches to \( \sigma_\varepsilon \), restarting the history from stage \( t_\varepsilon \). Formally:
\begin{itemize}
    \item For \( 1 \leq t < t_\varepsilon \), \( \sigma_t(z_1, \ldots, z_t) = \sigma_\ast(z_1, \ldots, z_t) \).
    \item For \( t \geq t_\varepsilon \), \( \sigma_t(z_1, \ldots, z_t) = \sigma_\varepsilon(z_{t_\varepsilon}, \ldots, z_t) \).
\end{itemize}

Since \eqref{choosing-te} holds, choose \(\delta>0\) sufficiently small so that, for every \(i\in I\),
\begin{equation}\label{conditionondelta}
\|g^i\|\sum_{t=t_\varepsilon}^{\infty} w^i_t
+
\|g^i\|(n_\varepsilon+t_\varepsilon)(w^i_{t_\varepsilon}+\delta)
<\varepsilon_1,
\qquad
w^i_{t_\varepsilon}+\delta<\theta_i .
\end{equation}
We prove that \(\sigma\) is a robust \(\varepsilon\)-equilibrium for
\(B(w,\delta)\cap W^I\).
If \( u \in B(w, \delta) \cap W^I \) is an impatient valuation profile,
analogously to inequality (\ref{eqn:maxmintau1new}) (in the proof of Lemma~\ref{lem:uniformmaxminatwinoverlineW}), we have
\begin{align}
E^z_{\s} \left[ \sum_{1 \leq t < t_\varepsilon} u^i_t g^i_t \right]
&\geq v^i_\ast[z] - \beta_i E^z_{\s} v^i[z_{t_\varepsilon}] - \varepsilon_1 - \varepsilon_2
\label{eqn:tau1new}
\end{align}
(where the term \(-\varepsilon_2\) on the right hand side can be dropped as we use equality~(\ref{eqn:sigmaFt3new}) in place of inequality~(\ref{eqn:maxminsigmaitauminusi})).
Analogously to inequality~(\ref{eqn:maxmintau2new}) (in the proof of Lemma~\ref{lem:uniformmaxminatwinoverlineW}),
\begin{equation}{\label{eqn:sigmaFt5}}
E_{\s}^z \left[\sum_{t=t_\varepsilon}^\infty u^i_t g^i_t \right]
\geq \beta_i\, E_{\s}^zv^i[z_{t_\varepsilon}] - \varepsilon_0 - 3\varepsilon_1.
\end{equation}

The tail estimate \eqref{eqn:sigmaFt5} is obtained by applying Lemma~\ref{lem:cesaro_patient} directly to the non-normalized tail weights
\[
\lambda_r:=u^i_{t_\varepsilon+r-1},\qquad r\ge1 .
\]
Indeed, \((\lambda_r)_r\) is nonincreasing, \(\sum_r\lambda_r\le1\), and, by \eqref{conditionondelta},
\[
\lambda_1=u^i_{t_\varepsilon}\le w^i_{t_\varepsilon}+\delta<\theta_i .
\]
Applying Lemma~\ref{lem:cesaro_patient} after conditioning on \(H_{t_\varepsilon}\), with
\[
M=\|g^i\|,\qquad N=n_\varepsilon+1,\qquad
\eta=\varepsilon_0,\qquad a=v^i[z_{t_\varepsilon}],
\]
and with tolerance \(\varepsilon_0+2\varepsilon_1\), gives the lower tail
estimate with
\[
\sum_{t\ge t_\varepsilon}u^i_t v^i[z_{t_\varepsilon}].
\]
Finally,
\[
\left|
\sum_{t=t_\varepsilon}^{\infty}u^i_t-\beta_i
\right|
\le
\sum_{t=t_\varepsilon}^{\infty}w^i_t+t_\varepsilon\delta,
\]
and \eqref{conditionondelta} converts this into the additional \(\varepsilon_1\) error. The upper tail estimate used below is obtained in the same way, using the upper inequality in \((2.a^\ast)\) and the upper part of Lemma~\ref{lem:cesaro_patient}.

Therefore, by adding inequalities~\eqref{eqn:sigmaFt5} and \eqref{eqn:tau1new}, and
the definition of \( \gamma^i_u({\s}) \), we have
\[
\gamma^i_u({\s})[z] = E^z_{\s} \sum_{t=1}^\infty u^i_t g^i_t \geq v^i_\ast[z] - \varepsilon_0- 4\varepsilon_1 - \varepsilon_2\geq v^i_\ast[z] - \varepsilon/2.
\]

Fix a player \( i \) and a strategy \( \tau^i \) of player~\( i \).
Repeating the preceding argument with \(E^z_{\sigma}\) replaced by
\(E^z_{\sigma^{-i},\tau^i}\), and using the upper estimate
\eqref{eqn:sigmaFt3bnew} in place of \eqref{eqn:sigmaFt3new}, yields
\[
\gamma^i_u({\s}^{-i}, \tau^i)[z] = E^z_{{\s}^{-i}, \tau^i} \sum_{t=1}^\infty u^i_t g^i_t \leq v^i_\ast[z] + \varepsilon/2.
\]

Therefore, inequality~\eqref{SGrobustepsilonequilibriumatw} holds for
\(U=B(w,\delta)\cap W^I\). Hence \(\sigma\) is a robust
\(\varepsilon\)-equilibrium for \(B(w,\delta)\cap W^I\).
This proves Lemma~\ref{lem:uniformequilibriumatwinoverlineW}.

In particular, \(\Gamma\) has robust equilibria at every \(w\in\overline{W}^I\).

The proofs of the implications \((2.a^\ast) \implies (2.b) \implies (2.c)\) are analogous to the proofs of \((1.a)\ \implies\ (1.b) \ \implies\ (1.c)\).

Assume that condition \((2.a^\ast)\) holds. By Lemma~\ref{lem:uniformequilibriumatwinoverlineW}, \(\Gamma\) has robust \(\varepsilon\)-equilibria at every \(w\in\overline W^I\). Suppose that condition \((2.b)\) fails. Then there exist \(\varepsilon>0\) and a sequence \((U_m)_m\) of nonempty subsets of \(W^I\), with \(\operatorname{diam}(U_m)\to 0\), such that no \(U_m\) admits a robust \(\varepsilon\)-equilibrium. By the product version of Lemma~\ref{lem:accumulationofUn}, there exists \(w\in\overline W^I\) such that every neighborhood of \(w\) contains some \(U_m\). By Lemma~\ref{lem:uniformequilibriumatwinoverlineW}, there are \(\delta>0\) and a strategy profile \(\sigma\) that is a robust \(\varepsilon\)-equilibrium for \(B(w,\delta)\cap W^I\).
For some \(m\), \(U_m\subset B(w,\delta)\cap W^I\), a contradiction.
Therefore \((2.b)\) holds.

This completes the proof that condition \((2.a^\ast)\) implies condition \((2.b)\).

\smallskip
Assume that condition \((2.b)\) holds, and fix \(\varepsilon>0\). Choose
\(r>0\) such that every set \(U\subseteq W^I\) with
\(\operatorname{diam}(U)<r\) admits a robust \(\varepsilon\)-equilibrium. Let \(0<\delta<r/4\).

Since \(\overline W^I\) is compact and \(W^I\) is dense in \(\overline W^I\), there exist finitely many points \(w(1),\ldots,w(k)\in W^I\) such that the balls \(B(w(j),\delta)\), \(j=1,\ldots,k\), cover \(W^I\). Define
\[
U_j:=B(w(j),\delta)\cap W^I.
\]
Then \(U_j^\delta\subseteq B(w(j),2\delta)\cap W^I\), whose diameter is at most \(4\delta<r\). Hence, by condition \((2.b)\), there exists a strategy profile \(\sigma_j\) that is a robust \(\varepsilon\)-equilibrium for \(U_j^\delta\).
Since the sets \(U_j\) cover \(W^I\), condition \((2.c)\) follows.

This completes the proof that condition \((2.b)\) implies condition \((2.c)\).

\smallskip

Assume that condition \((2.c)\) holds. Then, for every \(\varepsilon>0\), there exists a finite list \((U_j,\sigma_j)_{j\in J}\), where each \(\sigma_j\) is a robust \(\varepsilon\)-equilibrium for \(U^\delta_j\), and the sets \(U_j\) cover \(W^I\).

Let  \(U\) be any set of profiles of impatient valuations with \(\operatorname{diam}(U)<\delta\). If \(U\cap U_j\neq \emptyset\), then \(U\subseteq U^\delta_j\), and hence \(\sigma_j\) is a robust \(\varepsilon\)-equilibrium for \(U\).

Therefore, condition \((2.c)\) implies condition \((2.b)\).

\smallskip
Now assume condition \((2.b)\). Fix \(\varepsilon>0\). By condition \((2.b)\), there exists \(\delta>0\) such that every set \(U\subseteq W^I\) with \(\operatorname{diam}(U)<\delta\) admits a robust \(\varepsilon\)-equilibrium. Choose \(n_\varepsilon\) so large that \(1/n_\varepsilon<\delta\), and set
\[
U_\varepsilon:=\{w\in W^I: w^i_1\le 1/n_\varepsilon
\text{ for every } i\in I\}.
\]
Then \(\operatorname{diam}(U_\varepsilon)\le 1/n_\varepsilon<\delta\). Let \(\sigma\) be a robust \(\varepsilon\)-equilibrium for \(U_\varepsilon\), and set
\[
v_\varepsilon:=\gamma_{n_\varepsilon}(\sigma).
\]
For every \(n\ge n_\varepsilon\), both \(u_n\) and \(u_{n_\varepsilon}\) belong to \(U_\varepsilon\). Hence, for every player \(i\), every initial state \(z\), and every unilateral deviation \(\tau^i\),
\[
\gamma^i_n(\sigma^{-i},\tau^i)[z]
\le \gamma^i_{n_\varepsilon}(\sigma)[z]+\varepsilon
= v^i_\varepsilon[z]+\varepsilon .
\]
Moreover, applying the robust inequality with the deviation \(\tau^i=\sigma^i\) gives
\[
\gamma^i_n(\sigma)[z]
\ge \gamma^i_{n_\varepsilon}(\sigma)[z]-\varepsilon
= v^i_\varepsilon[z]-\varepsilon .
\]
Thus condition \((2.a^\ast)\) holds.

This completes the proof that \((2.b)\) implies \((2.a^\ast)\),
and thus establishes the equivalence of all stated conditions in Theorem~\ref{thm:robustequilibrium}.

\subsection{Proof of Theorem~\ref{thm:characterisation_mediated_equilibria}.}

For each of the three mediation classes considered in Theorem~4, ordinary
strategy profiles are embedded as degenerate mediated profiles: the mediator's signals may simply be ignored. Moreover, the finite-prefix/continuation construction
used in Lemma~\ref{lem:uniformequilibriumatwinoverlineW} can be represented within the same mediation class. In the extensive-form correlated case this is immediate, since signals may be sent during the play.
In the sunspot and correlated cases, the initial public signal, in the sunspot case, or the initial private signals, in the correlated case, can select a complete contingent plan before play begins. This plan specifies the prefix behavior and, for each possible state reached at stage \(t_\varepsilon\), the continuation behavior to be used from that state.
Since the state at stage \(t_\varepsilon\) is observed in the stochastic game, the players can then follow the branch of the pre-specified continuation plan corresponding to that state. Thus the constructed mediated profile is again a sunspot, correlated, or extensive-form correlated strategy profile, respectively. Therefore condition \((M.a^\ast)\) implies local robust \(\mathcal M\) equilibria at every \(w\in\overline W^I\).

The compactness argument used in the proof of \((2.a^\ast)\Rightarrow(2.b)\), and the finite-cover argument used in the proof of \((2.b)\Rightarrow(2.c)\), then give
\[
(M.a^\ast)\Rightarrow(M.b)\Rightarrow(M.c).
\]
Moreover, \((M.b)\) and \((M.b^\ast)\) are equivalent by definition.

We prove \((M.c)\Rightarrow(M.d)\). Fix \(\varepsilon>0\), and choose
\(0<\eta<\varepsilon\). Apply \((M.c)\) with tolerance \(\eta\). Thus there are \(\delta_0>0\) and a finite family \((\sigma_j,U_j)_{j\in J}\), where each \(\sigma_j\) is an \(\mathcal M\)-strategy profile and the sets \(U_j\) cover \(W^I\), such that \(\sigma_j\) is a robust \(\mathcal M\) \(\eta\)-equilibrium for \(U_j^{\delta_0}\).

For each \(j\), let
\[
V_j:=\{w\in\overline W^I:d(w,U_j)<\delta_0/2\}.
\]
The sets \(V_j\) cover \(\overline W^I\). Choose continuous functions
\(f_j:\overline W^I\to\mathbb R_+\) such that \(f_j\) is positive on \(V_j\) and vanishes outside \(V_j\), and set
\[
\alpha_j(w):=\frac{f_j(w)}{\sum_{\ell\in J}f_\ell(w)} .
\]
Then \((\alpha_j)_{j\in J}\) is a continuous partition of unity. For
\(w\in W^I\), define
\[
\sigma_w:=\sum_{j\in J}\alpha_j(w)\sigma_j .
\]
This convex combination is implemented by drawing an auxiliary label \(j\) before play begins, with probabilities \((\alpha_j(w))_{j\in J}\), and then using the mediated strategy profile \(\sigma_j\). Since \(\mathcal M\) is closed under mixing, \(\sigma_w\) is an \(\mathcal M\)-strategy profile.

If \(\alpha_j(w)>0\), then \(w\in V_j\). Hence, for every
\(u,u'\in B(w,\delta_0/2)\cap W^I\), we have \(u,u'\in U_j^{\delta_0}\). Therefore each active \(\sigma_j\) is a robust \(\mathcal M\) \(\eta\)-equilibrium for \(B(w,\delta_0/2)\cap W^I\). By linearity
of the payoff and incentive inequalities under mediated mixtures, \(\sigma_w\) is a robust \(\mathcal M\) \(\eta\)-equilibrium, hence a robust \(\mathcal M\) \(\varepsilon\)-equilibrium, for
\[
B(w,\delta_0/2)\cap W^I .
\]
Finally, by the metric convention for mediated strategies in Section~3.2, the map \(w\mapsto\sigma_w\) is uniformly continuous, because it is obtained from a fixed finite list \((\sigma_j)_{j\in J}\) with uniformly continuous weights \((\alpha_j(w))_{j\in J}\). This proves \((M.c)\Rightarrow(M.d)\).

It remains to prove \((M.d)\Rightarrow(M.a^\ast)\). Fix \(\varepsilon>0\), and apply \((M.d)\) with tolerance \(\varepsilon\). Let \(\delta>0\) and \(w\mapsto\sigma_w\) be as in \((M.d)\). Choose \(n_\varepsilon\) so large that \(1/n_\varepsilon<\delta\), and write
\[
{\bf u}_n:=(u_n,\ldots,u_n)\in W^I .
\]
Set
\[
\sigma_\varepsilon:=\sigma_{{\bf u}_{n_\varepsilon}},
\qquad
v_\varepsilon:=\gamma_{n_\varepsilon}(\sigma_\varepsilon).
\]
For every \(n\ge n_\varepsilon\), both \({\bf u}_n\) and
\({\bf u}_{n_\varepsilon}\) belong to
\(B({\bf u}_{n_\varepsilon},\delta)\cap W^I\). Hence, for every
\(k,m\ge n_\varepsilon\), every player \(i\), every initial state \(z\), and every admissible unilateral deviation \(\tau^i\),
\[
\gamma^i_k(\sigma_\varepsilon^{-i},\tau^i)[z]
\le
\gamma^i_m(\sigma_\varepsilon)[z]+\varepsilon .
\]

Taking \(k=n\) and \(m=n_\varepsilon\) gives, for every
\(n\ge n_\varepsilon\),
\[
\gamma^i_n(\sigma_\varepsilon^{-i},\tau^i)[z]
\le
\gamma^i_{n_\varepsilon}(\sigma_\varepsilon)[z]+\varepsilon
=
v^i_\varepsilon[z]+\varepsilon .
\]
This is the upper incentive inequality in \((M.a^\ast)\).

Taking \(k=n_\varepsilon\), \(m=n\), and \(\tau^i=\sigma^i_\varepsilon\), gives
\[
\gamma^i_{n_\varepsilon}(\sigma_\varepsilon)[z]
=
\gamma^i_{n_\varepsilon}(\sigma_\varepsilon^{-i},\sigma^i_\varepsilon)[z]
\le
\gamma^i_n(\sigma_\varepsilon)[z]+\varepsilon .
\]
Therefore,
\[
\gamma^i_n(\sigma_\varepsilon)[z]
\ge
\gamma^i_{n_\varepsilon}(\sigma_\varepsilon)[z]-\varepsilon
=
v^i_\varepsilon[z]-\varepsilon .
\]
This is the lower payoff inequality in \((M.a^\ast)\). Thus condition
\((M.a^\ast)\) holds. This completes the proof of Theorem~\ref{thm:characterisation_mediated_equilibria}.

\section{Robustness in Other Models of Multistage Games}\label{sec:othermodels}

The notions of uniform equilibrium payoff and robust equilibrium, developed in this paper for discrete-time (DT) stochastic games, extend naturally to broader classes of dynamic games.
These include the general model of repeated games (which is a unified framework that captures repeated interactions with incomplete information and imperfect monitoring; see \cite{mersorzam}), models with a continuum of states, and continuous-time (CT) stochastic games.

In what follows, we summarize a few results about the existence of robust optimization and robust equilibria in several such multistage game settings.

We begin with new results on repeated games with incomplete information, then proceed to several classes of stochastic games (in both discrete and continuous time), and conclude with examples that highlight the limits of robust optimization when no uniform value exists.

\subsection{Repeated games with incomplete information}
In a two-player zero-sum repeated game with incomplete information on one side, nature selects a state \(z\in S\), where \(S\) is a finite set of
states, according to a commonly known prior \(p\in \Delta(S)\). Player~1 observes \(z\), while Player~2 remains uninformed.

At each stage \(t\), Player~1 selects an action \(a_t\in A\), and Player~2 simultaneously selects \(b_t\in B\), where \(A\) and \(B\) are finite action sets.

After each stage, both players observe the action taken by the other player. The stage payoff to Player 1 at stage $t$ is \(g(z,a_t,b_t)\), where \(g:S\times A\times B\to \mathbb{R}\) is a given payoff function.

A strategy profile \((\s,\tau)\), together with the prior \(p\), induces a probability distribution \(P_{\s,\tau}^p\) over the infinite histories  \((z,a_1,b_1,a_2,b_2,\ldots)\).  Expectation with respect to  this distribution is denoted by \(E_{\s,\tau}^p\).

Given an impatient valuation \(w\in W\), define the game \(\Gamma_w(p)\), in which Player~1 aims to maximize the weighted average payoff \[\gamma_w(\s,\tau)[p]=E_{\s,\tau}^p \sum_{t=1}^\infty w_t \cdot g(z,a_t,b_t).\]

The game \(\Gamma_w(p)\) has a value, \(v(w)[p]\), and both players have optimal strategies.

\begin{defi}Let \(w\in \overline{W}\). We say that the game \(\Gamma(p)\) has \emph{robust optimal strategies at \(w\)} if
for every \( \ve>0\), there exist \(\delta>0\) and strategies \(\s_\ve,\tau_\ve\) such that for all \(\s, \tau\), and all \(u,u'\in W\cap B(w,\delta)\), we have
\begin{align} \gamma_u(\s_\ve,\tau)[p]& \geq  v(u')[p]-\ve, \label{SGd}\
\gamma_u(\s,\tau_\ve)[p] &\leq    v(u')[p] +\ve.
\end{align}
\end{defi}

It is known that these games admit a uniform value. However, due to the asymmetry of information---specifically, Player 2 not observing the state---the setting does not fall within the scope of Theorem \ref{thm:robustoptimalstrategy}. In particular, one cannot conclude the existence of a robust  \(\ve\)-optimal strategy at boundary points \(w\in \overline{W}\setminus W\) using prior results. We record the following results without proof:

\begin{result}Every two-player zero-sum repeated game with incomplete information on one side has robust optimal strategies.
\end{result}

\begin{result}Every two-player repeated game with incomplete information on one side has robust equilibria.
\end{result}

The notion of \emph{robust equilibria} in this setting is defined analogously to its counterpart in stochastic games. Specifically, the game \(\Gamma_p\) is said to have robust equilibria if, for every \(\varepsilon>0\), there exists $\delta>0$ such that for any set $U$ of profiles of impatient valuations with a diameter of less than $\delta$, there exists a strategy profile $\sigma$ that is an
\(\varepsilon\)-equilibrium in
\(\Gamma_p(u)\)  for every \(u\in U\), and such that the variation in the resulting payoffs across all $u\in U$ is at most $\varepsilon$.

\subsection{Two-person non-zero-sum repeated games with symmetric incomplete information}
\cite{neysor98} proved  that these games admit a uniform equilibrium payoff. The analysis of such games reduces to a stochastic game with a continuum of states, where the state space corresponds to players'  posterior beliefs. Since Theorem \ref{thm:robustequilibrium} holds also for stochastic games with a continuum state space, it follows that these games admit robust equilibria.

\subsection{CT non-zero-sum finite stochastic games}
An impatient valuation of a continuous stream of payoffs $[0,\infty)\ni t\mapsto g_t$, is an integral $\int_0^\infty w_tg_t\,dt$, where  \(w:[0,\infty)\to \mathbb{R}_+\) is a nonincreasing function with \(\int_0^\infty w_t\,dt= 1\). Let \(W_{ct}\) and \(\overline{W}_{ct}\) denote the spaces of all nonincreasing densities \(w:[0,\infty)\to \mathbb{R}_+\) with \(\int_0^\infty w_t\,dt= 1\) and  \(\int_0^\infty w_t\,dt\leq 1\), respectively. We endow these spaces with the metric
\[d(u,w):=\operatorname{ess\,sup}_{t\geq 0}|u_t-w_t|.\]

A strategy profile  \(\s\) in a CT stochastic game \(\Gamma\) is a robust \(\ve\)-equilibrium for a set  \(U\subseteq W^I_{ct}\) of profiles of impatient valuations if, for every player \(i\),  every strategy \(\tau^i\) of player \(i\), and all \(u,w\in U\), we have:
\[E_{\s^{-i},\tau^i}\int_0^\infty u^i_tg^i_t\,dt\leq E_\s\int_0^\infty w^i_tg^i_t\,dt+\ve.\]
We say that the game \(\Gamma\) \emph{admits robust equilibria} if for every \(\varepsilon>0\), there exists \(\delta>0\) such that for every set \(U\subseteq W^I_{ct}\) with a diameter of less than \(\delta\), there exists a robust \(\ve\)-equilibrium for \(U\).

\cite{neyman2017} proves that every CT stochastic game with finitely many states and actions admits a uniform equilibrium payoff. Building on this result---and using the techniques developed in this paper, along with tools such as jointly controlled lotteries---one can show that every such game admits robust equilibria. Moreover, for every \(\varepsilon>0\), there exists a uniformly continuous map  \(w\mapsto \sigma_{w,\varepsilon}\) from the space of impatient valuations \(W^I_{ct}\) to  the space of strategy profiles, such that \( \sigma_{w,\varepsilon}\) is an \(\varepsilon\)-equilibrium of \(\Gamma_w\).

\subsection{Zero-sum absorbing games with compact action spaces}
These are stochastic games with a single nonabsorbing state, where a state \(z\) is absorbing if \(q(z\mid z,a)=1\) \(\forall a\in A\).
\cite{merneyros} prove that a zero-sum absorbing game with compact action spaces has a uniform value; hence, by Theorem \ref{thm:robustoptimalstrategy}, it has robust optimal strategies.

\subsection{Limits of Robust Optimization: Games without a Uniform Value}
In two-player zero-sum repeated games, the absence of a uniform value precludes the existence of robust optimal strategies. We briefly describe several settings in which no uniform value exists and hence robust optimization cannot be guaranteed.

\begin{itemize}
\item \citet{vig13} constructs a zero-sum stochastic game with compact action spaces in which the asymptotic value fails to exist; in particular, the uniform value does not exist.

\item \citet{zil} proves that in stochastic games with finitely many states and actions, observable actions, and unobservable states, the asymptotic value may fail to exist.

\item \citet{aummas95} show that two-player zero-sum repeated games with incomplete information on both sides need not have a uniform value.

\item One can easily construct stochastic games with countably infinite state spaces and finitely many actions in which the asymptotic value does not exist.
\end{itemize}

\section{Robustness to imprecision in the specification of the transition probabilities of finite stochastic games}\label{sec:transitionsrobustness}

We conclude by highlighting an open question concerning robustness to small imprecisions in the transition probabilities of a finite stochastic game.

The strategy space and the space of play in a two-player zero-sum finite stochastic game \(\langle I, S, A, g, q\rangle\) are determined by the sets \(S\) and \(A\), and are independent of the payoff function \(g\) and the transition function \(q\). Let \(\sigma, \sigma_\varepsilon\) denote strategies of Player~1, and \(\tau, \tau_\varepsilon\) denote strategies of Player~2. Fix \(\Gamma^q = \langle I, S, A, g, q\rangle\), and consider a sequence of games \(\Gamma_k = \langle I, S, A, g, q_k \rangle\) such that:
\begin{itemize}
    \item \(q_k(z' \mid z, a) = 0\) whenever \(q(z' \mid z, a) = 0\), and
    \item \(q_k(z' \mid z, a) \to q(z' \mid z, a)\) whenever \(q(z' \mid z, a) > 0\).
\end{itemize}

\cite{solan03} proves that the sequence of uniform values of \(\Gamma_k\) converges to the uniform value \(v\) of \(\Gamma\). That is, for every \(\varepsilon > 0\), there exists \(k_\varepsilon\) such that for every \(k > k_\varepsilon\), there exist \((\sigma_\varepsilon, \tau_\varepsilon, n_\varepsilon)\) such that for all \(z \in S\), \(n > n_\varepsilon\), and strategy profiles \((\sigma, \tau)\), the following holds:
\[
\varepsilon + E^{k,z}_{\sigma_\varepsilon, \tau}[u_n(g)] \geq v[z] \geq E^{k,z}_{\sigma, \tau_\varepsilon}[u_n(g)]-\varepsilon,
\]
where \(E^{k,z}_{\sigma, \tau}\) denotes the expectation under the law of play induced by \((\sigma, \tau)\) in the game \(\Gamma_k\) starting at state \(z\).

In this result, the strategy profiles \((\sigma_\varepsilon, \tau_\varepsilon)\) and the horizon \(n_\varepsilon\) may depend on \(k\). The following open question asks whether this dependence can be removed.

\paragraph{\bf Open Problem.} Does there exist a functions \(v_w: S \to \mathbb{R}\), \(w\in W\), such that for every \(\varepsilon > 0\), there exist \((\sigma_\varepsilon, \tau_\varepsilon, n_\varepsilon, k_\varepsilon)\) such that for all \(z \in S\), \(n > n_\varepsilon\), \(k > k_\varepsilon\), and strategies \((\sigma, \tau)\), we have
\[
\varepsilon + E^{k,z}_{\sigma_\varepsilon, \tau}[u_n(g)] \geq v_w[z] \geq E^{k,z}_{\sigma, \tau_\varepsilon}[u_n(g)]-\varepsilon.
\]

A positive answer would imply the following strong form of robustness: for every \(w \in \overline{W}\) and \(\varepsilon > 0\), there exist \((\sigma_\varepsilon, \tau_\varepsilon, n_\varepsilon, k_\varepsilon)\) such that for all \(k, k' > k_\varepsilon\), \(n > n_\varepsilon\), \(z \in S\), \(w \in \overline{W}\), and impatient valuations \(u, u' \in B(w, 1/k_\varepsilon) \cap W\), we have:
\[
\varepsilon + E^{k,z}_{\sigma_\varepsilon, \tau}[u(g)] \geq v_w[z] \geq E^{k',z}_{\sigma, \tau_\varepsilon}[u'(g)]-\varepsilon.
\]

This type of robustness is theoretically significant because real-world models often involve imprecise transition probabilities, e.g. due to estimation errors, evolving environments, or misspecification. The ability to select approximately optimal strategies that remain effective across nearby transition structures would yield strong guarantees for decision-making under model uncertainty. It would also advance our understanding of the structural stability of strategic behavior in dynamic environments.

\bibliographystyle{plainnat}%
\bibliography{references}
\appendix
\section{Symmetric Markov equilibria of stochastic games with weighted average payoffs}\label{appendixMarkov}
This appendix establishes the existence of symmetric Markov equilibria in finite stochastic games with weighted average payoffs. As a corollary, it implies that every finite stochastic game with a profile of impatient valuations admits a Markov equilibrium. This result, which is not used directly in our arguments, highlights the contrast with the necessity of non-Markovian strategies for achieving robust optimization.

The proof presented here follows standard arguments in the theory of dynamic games. However, as we are unaware of a direct reference for this result in the existing literature, we provide a self-contained derivation for completeness.

For each $s\in \mathbb{N}$ and a profile $w=(w^i)_{i\in I}$ of averaging weights, i.e., $w^i=(w^i_1,w^i_2,\ldots)$ with $w^i_t\geq 0$ and $\sum_{t=1}^\infty w^i_t=1$, we set  $^s\!w$ to be the profile  $(^s\!w^i)_{i\in I}$ of averaging weights, where $^s\!w^i_t=w^i_{s+t-1}/\sum_{t\geq s}w^i_t$ if $\sum_{t\geq s}w^i_t>0$ and $^s\!w^i_t=1_{\{t=1\}}$ if $\sum_{t\geq s}w^i_t=0$. Let $w=(w^i)_{i\in I}$ be a profile of averaging weights and let $W_w$ be the set $\{^s\!w: s\in \mathbb{N}\}$. Note that each element of $W_w$ is a profile of averaging weights.

Let $\Gamma$ be a non-zero-sum stochastic game with finitely many states and actions with a set of players $I$, and let $w$ be a profile of averaging weights. Let $\Gamma_w$ denote the stochastic game with the same state and action spaces as $\Gamma$, where the payoff to player $i\in I$ under the strategy profile $\s$ is $\gamma^i(\s):=E_\s \sum_{t=1}^\infty w^i_tg^i_t$.
Note that if $^s\!w =\,^tw$ then the subgames following the  plays $(z_1,\ldots,z_s=z)$ and $(z'_1,\ldots,z_t=z)$ are isomorphic.

A {\em Markov strategy} of $\Gamma_w$ is a strategy $\s$ such that $\s(z_1,a_1,\ldots,z_s)$ is a function of just $z_s$ and $s$.
A {\em symmetric Markov strategy} of $\Gamma_w$ is a strategy $\s$ such that $\s(z_1,a_1,\ldots,z_s)$ is a function of just $z_s$ and $^s\!w$.

\begin{thm}For every profile $w$ of averaging weights and a non-zero-sum stochastic game $\Gamma$ with finitely many states and actions, $\Gamma_w$ has a subgame perfect equilibrium that consists of symmetric Markov strategies.
\end{thm}
\bproof
Set $X_{i,z,s}=\Delta(A^i(z))$, $i\in I, z\in S,s\geq 1$.%
Note that $X_{i,z,s}$%
is a convex compact subset of a Euclidean space. The $(i,z,s)$-coordinate of an element $x\in \widehat{X}:=\times_{i\in I,z\in S, s\geq 1}X_{i,z,s}$ is denoted by $x_{i,z,s}$. Note that an element $x\in \widehat{X}$ represents a profile of Markov strategies $\s$, where $\s^i(z_1,\ldots,z_s)=x_{i,z,s}$.

We endow $\widehat{X}$  with the product topology, i.e., the minimal topology under which the projection from $\widehat{X}$ on its $(i,z,s)$-coordinate is continuous. The topological space $\widehat{X}$ is a (countably infinite) product of convex compact subsets of Euclidean spaces. Hence, $\widehat{X}$ is a convex and compact (and metrizable) subset of a locally convex topological vector space.

Let $X$ denote the subset of $\widehat{X}$ consisting of all strategy profiles $x\in \widehat{X}$ whose $(i,z,s)$-coordinate $x_{i,z,s}$ coincides with its $(i,z,t)$-coordinate $x_{i,z,t}$ whenever $^s\!w =\, ^tw$. Note that an element $x\in X$ corresponds to the profile $\s$ of symmetric Markov strategies, where $\s^i(z_1,\ldots,z_s)$ is the $(i,z,s)$-coordinate of $x$.

Recall that $\|g\|=\max_{i\in I,z\in S,a\in A}|g^i(z,a)|$.%

Let $Y_{i,z,s}=[-\|g\|,\|g\|]$, $\widehat{Y}:=\times_{i\in I,z\in S, s\geq 1}Y_{i,z,s}$. The $(i,z,s)$-coordinate of an element $y\in \widehat{Y}$ is denoted by $y_{i,z,s}$. We endow $\widehat{Y}$ with the product topology. The topological space $\widehat{Y}$ is a (countably infinite) product of convex compact subsets of Euclidean spaces. Hence, $\widehat{Y}$ is a convex and compact (and metrizable) subset of a locally convex topological vector space.

Let $Y$ denote the subset of $\widehat{Y}$ consisting of all elements in $\widehat{Y}$ whose $(i,z,s)$-coordinate $y_{i,z,s}$ coincides with its $(i,z,t)$-coordinate $y_{i,z,t}$ whenever $^s\!w=\,^t\!w$. Note that $Y$ is a convex closed subset $\widehat{Y}$, and hence $Y$ is a convex and compact (and metrizable) subset of a locally convex topological vector space.

We define a set-valued mapping $f=(f^X,f^Y)$ from $X\times Y$ to $X\times Y$ by defining a single-valued function $f^Y$ from $X\times Y$ to $Y$ and a set-valued function $f^X$ from $X\times Y$ to $X$ and by setting $f(x,y)=f^X(x,y)\times \{f^Y(x,y)\}$. Formally:

\[f^Y_{i,z,s}(x,y)=\max_{x^i\in X_{i,z,s}}(^s\!w^i_1g^i(z,x^{-i}_{z,s},x^i)+(1-\,^s\!w^i_1)\sum_{z'\in S}q(z,x^{-i}_{z,s},x^i)(z')y_{i,z',s+1}),\]

where  $x^{-i}_{z,s}$ stands for the list of $x_{i',z,s}$ with $i'\neq i$, and
\[f^X_{i,z,s}(x,y)=\arg\!\max_{x^i\in X_{i,z,s}}(^s\!w^i_1g^i(z,x^{-i}_{z,s},x^i)+(1-^s\!w^i_1)\sum_{z'\in S}q(z,x^{-i}_{z,s},x^i)(z')y_{i,z',s+1}).\]

The product space $X\times Y$ is a convex compact subset of a locally convex topological vector space, $f(X\times Y)\subseteq X\times Y$, and $f$ is upper semicontinuous and convex-valued. Therefore, by the Glicksberg fixed point theorem, $f$ has a fixed point $(x,y)$. As $x\in X$, $x$ is a profile of symmetric Markov strategies.

We now need to verify that $\sigma=x$ is a subgame perfect equilibrium strategy of $\Gamma_w$.
First, note that as $(x,y)$ is a fixed point of $f$,
\[E_\s (w^i_sg^i_s+\sum_{t>s}w^i_t y_{i,z_{s+1},s+1}\mid z_1,\ldots,z_s)=\sum_{t\geq s}w^i_t y_{i,z_s,s},\]
and therefore, for every $s'\leq s$, we have
\[E_\s (w^i_sg^i_s\mid z_1,\ldots,z_{s'})=E_\s (\sum_{t\geq s}w^i_t y_{i,z_s,s}-\sum_{t>s}w^i_t y_{i,z_{s+1},s+1}\mid z_1,\ldots,z_{s'}).\]
Summing these equalities over $s\geq s'$, we have
\[E_\s \sum_{s\geq s'}(w^i_sg^i_s\mid z_1,\ldots,z_{s'})=\sum_{s\geq s'}w^i_s\, y_{i,z_{s'},s'}.\]
In the last summation argument we used the fact that  $\sum_{s=s'}^\infty (\sum_{t\geq s}w^i_t y_{i,z_s,s}-\sum_{t>s}w^i_t y_{i,z_{s+1},s+1})$ is a telescopic sum and $\sum_{t>s}w^i_t y_{i,z_{s+1},s+1}$ converges to 0 as $s$ goes to infinity. Therefore, $\sum_{s=s'}^\infty (\sum_{t\geq s}w^i_t y_{i,z_s,s}-\sum_{t>s}w^i_t y_{i,z_{s+1},s+1})=\sum_{t=s'}^\infty w^i_t y_{i,z_{s'},s'}$.

Similarly, if $\tau^i$ is a strategy of player $i$ we have
\[E_{\tau^i,\s^{-i}} (w^i_sg^i_s+\sum_{t>s}w^i_t y_{i,z_{s+1},s+1}\mid z_1,\ldots,z_s)\leq \sum_{t\geq s}w^i_t y_{i,z_s,s},\]
and therefore, for every $s'\leq s$, we have
\[E_{\tau^i,\s^{-i}} (w^i_sg^i_s\mid z_1,\ldots,z_{s'})\leq E_{\tau^i,\s^{-i}}(\sum_{t\geq s}w^i_t y_{i,z_s,s}-\sum_{t>s}w^i_t y_{i,z_{s+1},s+1}\mid z_1,\ldots,z_{s'}).\]
Summing these inequalities over $s\geq s'$, we have
\[E_{\tau^i,\s^{-i}} \sum_{s\geq s'}(w^i_sg^i_s\mid z_1,\ldots,z_{s'})\leq \sum_{s\geq s'}w^i_s y_{i,z_{s'},s'}.\]
This shows that $\s=x$ is a subgame perfect equilibrium. As $x\in X$, $\s$ is a profile of symmetric Markov strategies.
\eproof

This result shows that, for arbitrary profiles of fixed averaging weights, symmetric Markov strategies suffice for equilibrium. This is in sharp contrast to the non-Markovian strategies typically required for robust optimization under small imprecision with respect to the players' averaging weights.

\end{document}